\documentclass[12pt]{article}

\usepackage{amssymb}
\usepackage{amsmath}
\usepackage{amscd}
\usepackage[dvips]{graphicx}
\usepackage{amsfonts,amsthm}

\numberwithin{equation}{section}

\newtheorem{theorem}{Theorem}[section]

\newtheorem{proposition}{Proposition}

\begin{document}

\title{{\large ON A PROCESSOR SHARING QUEUE\\ THAT MODELS BALKING}}

\author{
{\small Qiang Zhen}\thanks{
Department of Mathematics, Statistics, and Computer Science,
University of Illinois at Chicago, 851 South Morgan (M/C 249),
Chicago, IL 60607-7045, USA.
 Email: \texttt{qzhen2@uic.edu}
}
\and
{\small Johan S. H. van Leeuwaarden}\thanks{
Department of Mathematics and Computer Science,
Eindhoven University of Technology,
Room HG 9.13, P.O. Box 513, 5600 MB Eindhoven, The Netherlands.
Email: \texttt{j.s.h.v.leeuwaarden@tue.nl}
}
\and
{\small Charles Knessl}\thanks{
Department of Mathematics, Statistics, and Computer Science,
University of Illinois at Chicago, 851 South Morgan (M/C 249),
Chicago, IL 60607-7045, USA.
Email: \texttt{knessl@uic.edu}\
\newline\indent\indent{\bf Acknowledgements:} Knessl was partly supported by NSF grant DMS 05-03745 and NSA grant H 98230-08-1-0102. Van Leeuwaarden
was supported by a VENI grant from The Netherlands Organization for Scientific Research
(NWO).
}}
\date{{\small March 09, 2009}}
\maketitle

\begin{abstract}
\noindent We consider the processor sharing $M/M/1$-PS queue which also models balking. A customer that arrives and sees $n$ others in the system ``balks" (i.e., decides not to enter) with probability $1-b_n$. If $b_n$ is inversely proportional to $n+1$, we obtain explicit expressions for a tagged customer's sojourn time distribution. We consider both the conditional distribution, conditioned on the number of other customers present when the tagged customer arrives, as well as the unconditional distribution. We then evaluate the results in various asymptotic limits. These include large time (tail behavior) and/or large $n$, lightly loaded systems where the arrival rate $\lambda\to 0$, and heavily loaded systems where $\lambda\to\infty$. We find that the asymptotic structure for the problem with balking is much different from the standard $M/M/1$-PS queue. We also discuss a perturbation method for deriving the asymptotics, which should apply to more general balking functions.
\end{abstract}

\section{Introduction}
Balking refers to the phenomenon that customers, when forced to wait for service, refuse to join the queue. First introduced by Haight \cite{haight1957}, balking can be specified by a probability distribution associated with the system state. Specifically, a customer that finds $n$ customers in the systems upon arrival, balks with probability $1-b_n$ and joins the queue with probability $b_n$. Haight considered several examples of balking functions, including $b_n=1/(n+1)$, $b_n=1\{n\leq K\}$ (with $1\{\cdot\}$ the indicator function), and $b_n=e^{-c n}$. The latter case was also studied in detail by Morse \cite{morse1958}. In this paper we shall investigate the effect of balking upon systems with processor sharing (PS).

A model for the round-robin scheduling mechanisms in time-shared computer systems, processor sharing was first introduced by Kleinrock \cite{KLanalysis}, and refers to the service discipline under which every customer gets a fair share of the server.
It is by now well known that PS is intimately related to random order of service (ROS), which refers to the discipline where customers are chosen for service at random. First studied by Vaulot \cite{VA} and Pollaczek \cite{PO}, the ROS discipline has a long tradition in queueing theory. Pollaczek obtained the Laplace transform of the distribution of the steady-state waiting time $\mathcal{W}$ in the $M/M/1$-ROS queue, by solving a differential-difference equation. In fact, the latter was almost identical to the differential-difference equation studied by Coffman, Muntz and Trotter \cite{COF} for the $M/M/1$-PS queue. Indeed, by comparing these differential-difference equations, it is readily established that (see Cohen \cite{COH})
\begin{equation}
{\rm Prob}[\mathcal{V}>t]=C \cdot{\rm Prob}[\mathcal{W}>t]
\end{equation}
with $\mathcal{V}$ the steady-state sojourn time in the $M/M/1$-PS queue and  $C$ a constant. A probabilistic argument based on coupling was given in Borst, Boxma, Morrison and N{\'u}{\~n}ez-Queija \cite{BO}, and the equivalence result was shown to extend to other models as well, including finite capacity queues, repairman problems and networks. We shall show that for the $M/M/1$ queue with balking the equivalence result also holds.

The distribution of $\mathcal{W}$ does not have a simple representation. Pollaczek \cite{PO} was able to invert the Laplace transform and obtained a rather intricate but explicit integral representation for the waiting time distribution. The integral, along with the method of steepest descent, allowed Pollaczek to derive an intriguing asymptotic expression for the tail distribution ${\rm Prob}[\mathcal{W}>t]$. This asymptotic expression was rediscovered by Flatto \cite{FL}. Morrison \cite{MO} considered the heavy-traffic limit, where the traffic intensity $\rho\rightarrow 1$, and derived asymptotic expansions in powers of $1-\rho$ for the sojourn time distributions. The tail results of Pollaczek and Flatto were related to the heavy-traffic results in Morrison recently in Zhen and Knessl \cite{ZH}.

In this paper we consider the $M/M/1$-PS (or ROS) queue with balking. The sojourn time distribution (or waiting time distribution in the ROS model) satisfies a differential-difference equation that differs only slightly from the one considered by Pollaczek for the systems without balking. However, the analysis, and also the system behavior, changes drastically. We shall assume that $b_n=1/(n+1)$, so that the non-balking probability exactly matches the  share of a server that the customers gets upon arrival. For this choice of $b_n$ the differential-difference equation allows for an exact and asymptotic analysis. We obtain the following results:
\begin{itemize}
\item[(i)] An exact spectral representation for the sojourn time density in terms of  generalized Laguerre polynomials.
\item[(ii)] An expression for the Laplace-Stieltjes transform of the sojourn time density.
\item[(iii)]  Asymptotic results for tail probabilities when $\rho$ is fixed; asymptotics for the light-traffic case where $\rho\rightarrow 0$; and asymptotics for the heavy-traffic case where $\rho\rightarrow\infty$.
\end{itemize}
The heavy-traffic asymptotics are derived using a singular perturbation approach. It is also explained how this approach might be useful for analyzing models with more general balking functions.

\subsection{Equivalence relation}

We denote the sojourn time of a
non-balking customer that arrives to a PS queue with $n$ other customers
competing for service by $\mathcal{V}_n$, and the waiting
time of a non-balking customer that arrives to a ROS queue
with $n$ other customers waiting for service and one
additional customer in service by $\mathcal{W}_n$. Then we let $b_n$ and $b_n^r$ be the non-balking probabilities in the PS system and the ROS system, respectively, when there are $n$ customers in the system (including the customer in service).

\begin{proposition} If $b_0^r=1$ and $b_n^r=b_{n-1}$, $n=1,2,\ldots$,
then
\begin{equation}\label{eqr34}
\mathcal{V}_n\stackrel{d}{=}\mathcal{W}_n
\end{equation} and
\begin{equation}\label{eqrel}
{\rm Prob}[\mathcal{V}>t]=C \cdot{\rm Prob}[\mathcal{W}>t]
\end{equation}
with
\begin{equation}
C=\frac{1}{\rho} \cdot \frac{1+\sum_{n=1}^{\infty}\rho^n b_0\cdots b_{n-1}}{1+\sum_{n=1}^{\infty}\rho^n b_0\cdots b_{n}}.
\end{equation}
\end{proposition}
\begin{proof}
Borst, Boxma, Morrison and N\'{u}\~{n}ez-Queija \cite{BO} made the
observation that whenever a service
completion occurs in the PS system, each of the customers
present is equally likely to be the one that
departs due to the memoryless property of the
exponential distribution. In that respect, the pool of
customers competing for service under PS behaves
exactly as the pool of customers waiting for service
under ROS. Note that the arrival processes in both systems can be coupled due to the assumption that
$b_0^r=1$ and $b_n^r=b_{n-1}$, $n=1,2,\ldots$. A similar coupling argument as in
\cite{BO} then yields \eqref{eqr34}.

Let $\mathcal{N}_{p}$ and $\mathcal{N}_{r}$
denote the number of customers  at arrival epochs in the PS system and ROS system, respectively.
Then, ${\rm Prob}[\mathcal{N}_{p}=n]=\pi_0\rho^n b_0\cdots b_{n-1}$ and
\begin{align}\label{sdfgre1}
{\rm Prob}[\mathcal{V}>t]&=\frac{\sum_{n=0}^{\infty}\rho^n b_0\cdots b_{n}{\rm Prob}[\mathcal{V}_{n}>t]}{\sum_{n=0}^{\infty}\rho^n b_0\cdots b_{n}}.
\end{align}
Similarly,
${\rm Prob}[\mathcal{N}_{r}=n]=\pi_0^r\rho^n b_0^r\cdots b_{n-1}^r$ and
\begin{align}\label{sdfgre2}
{\rm Prob}[\mathcal{W}>t]&=\frac{\sum_{n=0}^{\infty}\rho^{n+1} b_0^r\cdots  b_{n+1}^r{\rm Prob}[\mathcal{W}_{n}>t]}{\sum_{n=0}^{\infty}\rho^{n} b_0^r\cdots  b_{n}^r}.
\end{align}
Upon comparing \eqref{sdfgre1} and \eqref{sdfgre2}, and using $b_0^r=1$, $b_n^r=b_{n-1}$ for $n=1,2,\ldots$,
and $\mathcal{V}_n\stackrel{d}{=}\mathcal{W}_n$, the equivalence relation \eqref{eqrel} follows.
\end{proof}

The $M/M/1/K$-PS queue can be viewed as a special case of the $M/M/1$-PS queue with balking by choosing $b_n=1$ if $n\leq K-1$ and 0 otherwise. In that case the equivalence relation becomes
\begin{equation}
{\rm Prob}[\mathcal{V}>t]=\frac{1}{\rho}\cdot \frac{1-\rho^{K+1}}{1-\rho^K} \cdot{\rm Prob}[\mathcal{W}>t],
\end{equation}
which was already obtained in \cite{BO}.
For this $M/M/1/K$-PS queue, Knessl \cite{knessl1993} uses singular perturbation techniques to construct asymptotic approximations
to the sojourn time distribution.

We shall consider the $M/M/1$-PS queue with balking $b_n=\frac{1}{n+1}$, in which case
\begin{equation}
{\rm Prob}[\mathcal{N}_p=n]=\frac{e^{-\rho}\rho^n}{n!}, \quad n=0,1,\ldots
\end{equation}
and
\begin{equation}\label{S1_equiv}
{\rm Prob}[\mathcal{V}>t]= \frac{e^\rho}{e^{\rho}-1}\cdot{\rm Prob}[\mathcal{W}>t].
\end{equation}

\section{Problem statement and summary of results}
We consider a processor sharing $M/M/1$ queue, which also models balking. Customers arrive at rate $\lambda$ and the service rate will be denoted by $\mu$. We can clearly scale time so as to have $\mu=1$, and then the traffic intensity is $\rho=\lambda/\mu=\lambda$. We let $\mathcal{V}_n$ be the sojourn time of a tagged customer that finds $n$ others in the system upon arrival. We then define
\begin{equation}\label{S2_cdf}
\mathbf{V}_n(t)=\mathrm{Prob}[\mathcal{V}_n>t].
\end{equation}
With the PS discipline, each customer receives service at rate $1/n$ when there are $n$ customers in service. 
When a tagged customer arrives we assume that he/she will enter the system with probability $b_{n}$, and ``balk" with the remaining probability $1-b_{n}$. 

The function $\mathbf{V}_n(t)$ satisfies the differential-difference equation
\begin{equation}\label{S2_recu1}
\mathbf{V}'_n(t)=\frac{n}{n+1}\,\mathbf{V}_{n-1}(t)-\Big(1+\rho\,b_{n}\Big)\,\mathbf{V}_n(t)+\rho\,b_{n}\mathbf{V}_{n+1}(t)
\end{equation}
with $\mathbf{V}_n(0)=1$. (There is a slight error in \cite{RI}, where this equation was previously given.) It is reasonable to define $b_0=1$ and have $b_n$ a decreasing function of $n$. Here we assume that
\begin{equation}\label{S2_bn}
b_n=\frac{1}{n+1},\;\;\;n=0,1,2,3,\cdots.
\end{equation}
Then we define the sojourn time density $p_n(t)$ by $p_n(t)=-\mathbf{V}'_n(t)$, which satisfies:
\begin{equation}\label{S2_recu2}
p'_n(t)=\frac{n}{n+1}\,p_{n-1}(t)-\Big(1+\frac{\rho}{n+1}\Big)\,p_n(t)+\frac{\rho}{n+1}\,p_{n+1}(t),\;\;\;t>0
\end{equation}
with the initial condition
\begin{equation}\label{S2_pn0}
p_n(0)=\frac{1}{n+1},\;\;n\geq 0.
\end{equation}
The above can be obtained by integrating (\ref{S2_recu1}) from $t=0$ to $t=\infty$ and using $\mathbf{V}_n(0)=1$.

Clearly (\ref{S2_bn}) is a very special case of $b_n$. We can consider also $b_n=\alpha/(n+1)$, since $\alpha$ can be incorporated into the traffic intensity $\rho$. However, with even a slight change (such as taking $b_n=1/(n+\beta)$ with $\beta\neq 1$) it seems that the problem is no longer amenable to exact solution. We shall discuss an asymptotic approach to solving (\ref{S2_recu1}) (cf. section 5), which should work also for more general $b_n$.

We give below various exact and asymptotic expressions for $p_n(t)$.

\begin{theorem} \label{th1}
The conditional sojourn time density has the following exact expression {\rm(}spectral representation{\rm)}:
\begin{equation}\label{S2_th1_p}
p_n(t)=\sum_{m=1}^\infty C_m(\nu_m)\,\phi_m(n,\nu_m)\,e^{\nu_m\,t}+\sum_{m=1}^\infty C_m(\widetilde{\nu}_m)\,\phi_m(n,\widetilde{\nu}_m)\,e^{\widetilde{\nu}_m\,t},
\end{equation}
where
\begin{equation}\label{S2_th1_nu}
\nu_m=-1+\frac{1}{2\,m}\Big[-\rho+\sqrt{\rho^2+4m\rho}\Big],
\end{equation}
\begin{equation}\label{S2_th1_nutilde}
\widetilde{\nu}_m=-1+\frac{1}{2\,m}\Big[-\rho-\sqrt{\rho^2+4m\rho}\Big],
\end{equation}
\begin{equation}\label{S2_th1_C}
C_m(\nu)=\frac{m^{m-1}}{m!}\,\frac{\nu}{\nu-1}\,e^{-m},
\end{equation}
and
\begin{equation}\label{S2_th1_phi}
\phi_m(n,\nu)=n!\,\Big(\frac{\nu+1}{-\rho}\Big)^n\,L^{(m-1-n)}_n\Big(\frac{\rho}{(\nu+1)^2}\Big).
\end{equation}
Here $L_n^{(\alpha)}(z)$ is the generalized Laguerre polynomial {\rm(}see {\rm \cite{MA})}.
\end{theorem}

If we take the Laplace transform of (\ref{S2_recu2}) and multiply by $n+1$, we have
\begin{equation}\label{S2_recu3}
\rho\,\widehat{p}_{n+1}(\theta)-\big[(n+1)(\theta+1)+\rho\big]\,\widehat{p}_n(\theta)+n\;\widehat{p}_{n-1}(\theta)=-1,
\end{equation}
where $\widehat{p}_n(\theta)=\int^\infty_0 p_n(t)e^{-\theta t}dt$. Solving the recurrence equation (\ref{S2_recu3}), we obtain another exact expression for $p_n(t)$, in terms of its Laplace transform.

\begin{theorem} \label{th2}
The Laplace-Stieltjes transform of the conditional sojourn time density has the following form:
\begin{equation}\label{S2_th2_phat}
\widehat{p}_n(\theta)=M\,G_n\sum_{l=0}^n\frac{\rho^l}{l!}\,H_l+M\,H_n\sum_{l=n+1}^\infty\frac{\rho^l}{l!}\,G_l,
\end{equation}
where
\begin{equation}\label{S2_th2_M}
M=M(\theta)\equiv \frac{\rho^{r+1}}{(1+\theta)\,\Gamma(r+1)}\,e^{r/\theta},
\end{equation}
\begin{equation}\label{S2_th2_G}
G_n=G_n(\theta)\equiv\int_0^{\frac{1}{(1+\theta)}}z^n\,\Big(\frac{1}{1+\theta}-z\Big)^r\exp\Big(-\frac{\rho\,z}{1+\theta}\Big) dz,
\end{equation}
\begin{equation}\label{S2_th2_H}
H_n=H_n(\theta)\equiv\int_{\frac{1}{(1+\theta)}}^\infty z^n\,\Big(z-\frac{1}{1+\theta}\Big)^r\exp\Big(-\frac{\rho\,z}{1+\theta}\Big) dz,
\end{equation}
and
$$r=r(\theta)\equiv\frac{\rho\,\theta}{(1+\theta)^2}.$$
The first two conditional moments of the sojourn time are
\begin{equation}\label{S2_th2_Mn}
\mathcal{M}_n=\int_0^\infty t\,p_n(t)\,dt=\frac{n+\rho}{2}+1,
\end{equation}
\begin{equation}\label{S2_th2_Sn}
\mathcal{S}_n=\int_0^\infty t^2\,p_n(t)\,dt=\frac{n^2}{3}+\Big(2+\frac{5}{6}\,\rho\Big)\,n+\frac{5}{6}\,\rho^2+3\rho+2.
\end{equation}

\end{theorem}

\vspace{2mm}

Using (\ref{S2_th2_phat}), we obtain the following asymptotic expansions for $p_n(t)$, valid for $\rho>0$ and $n$ and/or $t\to\infty$.
\begin{theorem} \label{th3}
For a fixed $\rho>0$ with $n,\;t\to\infty$, the conditional sojourn time density has the following asymptotic expansions:

\begin{enumerate}

\item $n\to\infty$, $n/t>1$,
\begin{equation}\label{S2_th3_1}
p_n(t)=\frac{1}{n}-\frac{\rho}{n(n-t)}+\frac{\rho-1}{n^2}+O(n^{-3}).
\end{equation}

\item $n/t=1+\Delta\,t^{-1/2}=1+O(t^{-1/2})$,
\begin{equation}\label{S2_th3_2}
p_n(t)\sim\frac{1}{2n}\,\mathrm{erfc}\Big(-\frac{\Delta}{\sqrt{2}}\Big)=\frac{1}{n\sqrt{\pi}}\int_{-\Delta/\sqrt{2}}^\infty \,e^{-u^2}\,du.
\end{equation}

\item $\Lambda_0<n/t<1$ with $\Lambda_0=\big(-\rho+\sqrt{\rho^2+4\rho}\;\big)/2$,
\begin{equation}\label{S2_th3_3}
p_n(t)\sim\frac{\Gamma(r_\ast+1)\,e^{r_\ast}}{\sqrt{2\pi}\,(1-n/t)^{r_\ast+1}}\,n^{-3/2-\rho\,t/n+\rho\,t^2/n^2}\,\Big(\frac{t}{n}\Big)^n\,e^{n-t},
\end{equation}
where
\begin{equation}\label{S2_th3_3_rstar}
r_\ast=r_\ast\Big(\frac{n}{t}\Big)\equiv\rho\,\frac{t^2}{n^2}\,\Big(\frac{n}{t}-1\Big).
\end{equation}

\item $n/t=\Lambda_0+\Lambda/\sqrt{t}$, $\Lambda=O(1)$,
\begin{equation}\label{S2_th3_4}
p_n(t) \sim \frac{\sqrt{\rho+4}-\sqrt{\rho}}{4\sqrt{\rho+4}}\,e^{-1}\,\Lambda_0^{-n}\,e^{-t+\Lambda_0\,t}\,\mathrm{erfc}\Big\{\frac{\Lambda}{\sqrt{2\Lambda_0}}\Big\}.
\end{equation}

\item $n/t<\Lambda_0$,
\begin{equation}\label{S2_th3_5}
p_n(t)\sim \frac{\sqrt{\rho+4}-\sqrt{\rho}}{2\sqrt{\rho+4}}\,e^{-1}\,\Lambda_0^{-n}\,e^{-t+\Lambda_0\,t}.
\end{equation}
\end{enumerate}
\end{theorem}

Expression (\ref{S2_th3_5}) applies also to $t\to\infty$ with $n=O(1)$, and gives the exponential decay rate of the density $p_n(t)$. We note that the right side of (\ref{S2_th3_5}) is precisely the $m=1$ term in the first sum in (\ref{S2_th1_p}), i.e., $C_1(\nu_1)\,\phi_1(n,\nu_1)\,e^{\nu_1\,t}$. If we start with a fixed large $n$ and increase time $t$ from $t=0$, we traverse cases 1-5 in Theorem 2.3 in the order given. The leading term in (\ref{S2_th3_1}) is $p_n(t)\sim 1/n$ for $t<n$ which corresponds to a uniform distribution. The $O(n^{-2})$ correction term(s) have a singularity as $t\uparrow n$, which indicates that the asymptotics become invalid. We also note that if $t=0$, (\ref{S2_th3_1}) becomes $p_n(0)=1/n-1/n^2+O(n^{-3})$, which is just the large $n$ expansion of the initial condition $p_n(0)=1/(n+1)$. As $t/n$ increases through one, there is a transition region (cf.\;(\ref{S2_th3_2})) and then for $t/n>1$ (but with $t/n<1/\Lambda_0$) the density becomes exponentially small, with a rather intricate dependence on the space-time ratio, as given in (\ref{S2_th3_3}). After another transition region where $t/n\approx 1/\Lambda_0$ (cf.\;(\ref{S2_th3_4})) the density becomes purely exponential in $t$, which corresponds to the dominant singularity in the Laplace transform $\widehat{p}_n(\theta)$, which occurs at $\theta=\nu_1<0$.

We next consider a small traffic intensity, $\rho\to 0^{+}$. We shall consider the time scales $t=O(\rho^{-1})$, $t=O(\rho^{-1/2})$ and $t=O(1)$.

\begin{theorem} \label{th4}
For $\rho\to 0^{+}$, the conditional sojourn time density has the following asymptotic expansions:

\begin{enumerate}
\item $t=\zeta/\rho=O(\rho^{-1})$,

\begin{enumerate}
\item $n=x/\rho=O(\rho^{-1})$ with $x>\zeta$,
\begin{equation}\label{S2_th4_11}
p_n(t)\sim 1/n.
\end{equation}

\item $n=x/\rho=O(\rho^{-1})$ with $x=\zeta+\Omega\,\sqrt{\rho}$, $\Omega=O(1)$
\begin{equation}\label{S2_th4_12}
p_n(t)\sim\frac{1}{2n}\,\mathrm{erfc}\Big(-\frac{\Omega}{\sqrt{2x}}\Big).
\end{equation}

\item $n=x/\rho=O(\rho^{-1})$ with $x<\zeta$,
\begin{equation}\label{S2_th4_13}
p_n(t)\sim\frac{\rho^{3/2}\,\zeta}{\sqrt{2\pi}\,(\zeta-x)}\,x^{-3/2}\,\exp\Big\{\frac{1}{\rho}\Big[x-\zeta+x\,\log(\zeta/x)\Big]\Big\}.
\end{equation}

\item $n=X/\sqrt{\rho}=O(\rho^{-1/2})$ with $X>\zeta$,
\begin{equation}\label{S2_th4_14}
p_n(t)\sim\frac{1}{\sqrt{2\pi}}\,\Gamma\Big(1-\frac{\zeta^2}{X^2}\Big)\,e^{-\zeta^2/X^2}\,n^{-3/2+\zeta^2/X^2}\,\Big(\frac{t}{n}\Big)^n\,e^{n-t}.
\end{equation}

\item $n=X/\sqrt{\rho}=O(\rho^{-1/2})$ with $X=\zeta+\rho^{1/4}\,Y$,
\begin{equation}\label{S2_th4_15}
p_n(t)\sim\frac{e^{-1}}{4}\,\rho^{-n/2}\,e^{-(1-\sqrt{\rho})\,t}\,\mathrm{erfc}\Big(\frac{Y}{\sqrt{2\zeta}}\Big).
\end{equation}

\item $n=X/\sqrt{\rho}=O(\rho^{-1/2})$ with $X<\zeta$,
\begin{equation}\label{S2_th4_16}
p_n(t)\sim\frac{e^{-1}}{2}\,\rho^{-n/2}\,e^{X/2}\,e^{-\zeta/2}\,e^{-(1-\sqrt{\rho})\,t}.
\end{equation}

\item $n=O(1)$,
\begin{equation}\label{S2_th4_17}
p_n(t)\sim\frac{e^{-1}}{2}\,\rho^{-n/2}\,e^{-\zeta/2}\,e^{-(1-\sqrt{\rho})\,t}.
\end{equation}
\end{enumerate}

\item $t=\omega/\sqrt{\rho}=O(\rho^{-1/2})$ and $n=O(1)$,
\begin{equation}\label{S2_th4_2}
p_n(t)\sim e^{-t}\,\rho^{-n/2}\,Q_n(\omega),
\end{equation}
where
\begin{eqnarray}\label{S2_th4_2Q}
Q_n(\omega)&=&\sum_{m=1}^\infty (-1)^n\,\frac{n!\,m^{m-n/2-1}}{2\,m!}\,e^{-m}\,L_n^{(m-1-n)}(m)\,e^{\omega/\sqrt{m}}\nonumber\\
&&+\sum_{m=1}^\infty \frac{n!\,m^{m-n/2-1}}{2\,m!}\,e^{-m}\,L_n^{(m-1-n)}(m)\,e^{-\omega/\sqrt{m}}.
\end{eqnarray}

\item $t,\;n=O(1)$,
\begin{equation}\label{S2_th4_3}
p_n(t)=p_n^{(0)}(t)+\rho\,p_n^{(1)}(t)+O(\rho^2),
\end{equation}
where
\begin{equation*}\label{S2_th4_30}
p_n^{(0)}(t)=\frac{1}{n+1}\,\frac{1}{2\pi i}\int_{Br}\Big[1-\Big(\frac{1}{1+\theta}\Big)^{n+1}\Big]\,\frac{e^{\theta\,t}}{\theta}\,dt=\frac{e^{-t}}{n+1}\,\sum_{l=0}^n\,\frac{t^l}{l!}
\end{equation*}
and
\begin{equation*}\label{S2_th4_31}
p_n^{(1)}(t)=\frac{e^{-t}}{n+1}\left[\frac{t^{n+2}}{(n+2)!}\,\sum_{l=0}^n \frac{1}{l+2}+\sum_{l=1}^{n+1}\frac{t^l}{l!}\Big(\frac{1}{n+2}-\frac{1}{n+2-l}\Big)\right].
\end{equation*}
Here $Br$ is a vertical Bromwich contour in the $\theta$-plane with $\Re(\theta)>0$.
\end{enumerate}
\end{theorem}

Some of the results in Theorem 2.4 for the time scale $t=O(\rho^{-1})$ can be derived as limiting cases of Theorem 2.3. However, this is not the case for the time ranges $t=O(\rho^{-1/2})$ and $t=O(1)$. For $n,\,t=O(1)$ the leading term in (\ref{S2_th4_3}) corresponds to the tagged customer entering the system and no further arrivals entering during his/her sojourn time.

The result in (\ref{S2_th4_2}) is obtained by letting $t=\omega/\sqrt{\rho}$ and taking $\rho\to 0$ in the exact expression (\ref{S2_th1_p}). If we let $\omega\to\infty$, the $m=1$ term in the first summation in (\ref{S2_th4_2Q}) dominates and this verifies the asymptotic matching between the scales $t=O(\rho^{-1})$ and $t=O(\rho^{-1/2})$, for $n=O(1)$.

Finally, we consider a large traffic intensity, $\rho\to\infty$. The structure of $p_n(t)$ is different in two cases.

\begin{theorem}\label{th5}
For $\rho\to\infty$, the conditional sojourn time density has the following asymptotic expansions:

\begin{enumerate}

\item $t=T\rho=O(\rho)$ and $n=N\rho=O(\rho)$,
\begin{equation}\label{S2_th5_1}
p_n(t)=\rho^{-1}\,P_0(N,T)+\rho^{-2}\,P_1(N,T)+O(\rho^{-3}),
\end{equation}
where
\begin{equation}\label{S2_th5_1p0}
P_0(N,T)=\frac{e^{U-T}}{N-U}=\frac{N-U-1}{(N-1)\,(N-U)}
\end{equation}
and $U=U(N,T)$ is defined implicitly by
\begin{equation}\label{S2_th5_1U}
\frac{U}{N-1}=1-e^{U-T}.
\end{equation}
If $N=1$ we obtain the explicit form $P_0(1,T)=e^{-T}$.

\item $t=\tau/\rho=O(\rho^{-1})$ and $n=O(1)$,
\begin{equation}\label{S2_th5_2}
p_n(t)\sim\int_0^1(1-\xi)^n\,J_0\Big(2\sqrt{\tau}\sqrt{-\xi-\log(1-\xi)}\;\Big)\,d\xi,
\end{equation}
where $J_0(\cdot)$ is the Bessel function of the first kind.

\end{enumerate}
\end{theorem}

We shall compute the correction term $P_1(N,T)$ in (\ref{S2_th5_1}) in section 5 and also give some alternate expressions for the leading term $P_0(N,T)$, as infinite series. The expression in (\ref{S2_th5_1}) remains valid for $n=O(1)$ and $t=O(\rho)$, as well as $t=O(1)$ and $n=O(\rho)$. For $N/T\gg 1$ we have $U\sim T$ and then $P_0(N,T)\sim 1/N$ which is consistent with $p_n(0)=1/(n+1)\sim \rho^{-1}/N$. For $T/N\gg 1$, $U\to -1$ and we obtain $P_0(N,T)\sim e^{-1}\,e^{-T}$, which is consistent with $C_1\,\phi_1\,e^{\nu_1\,t}$ for $\rho\to\infty$ and $t=O(\rho)$. Note that $\nu_1\sim -1/\rho$ from (\ref{S2_th1_nu}).

We remove the conditioning to get the unconditional sojourn time density for the PS model as
\begin{equation}\label{S2_p}
p_{_{PS}}(t)=\sum_{n=0}^\infty \frac{\rho^{n}}{n!}\,e^{-\rho}\,p_n(t).
\end{equation}
Then the density $p(t)$ for the ROS model follows from (\ref{S1_equiv}) as
$$p(t)=(1-e^{-\rho})\,p_{_{PS}}(t).$$
Note also that the full density, $p_{_{ROS}}(t)$, for the ROS model is $e^{-\rho}\,\delta(t)+p(t)$, since there is a non-zero probability that $\mathcal{W}=0$. The exact representation for $p_{_{PS}}(t)$ is as follows.

\begin{theorem}\label{th6}
The unconditional sojourn time density has the exact expression
\begin{equation}\label{S2_th_p}
p_{_{PS}}(t)=\sum_{m=1}^\infty C_m(\nu_m)\,\Phi_m(\nu_m)\,e^{\nu_m\,t}+\sum_{m=1}^\infty C_m(\widetilde{\nu}_m)\,\Phi_m(\widetilde{\nu}_m)\,e^{\widetilde{\nu}_m\,t},
\end{equation}
where
$$\Phi_m(\nu)=e^{-\rho}\,(-\nu)^{m-1}\,\exp\Big(\frac{\rho}{\nu+1}\Big).$$

\end{theorem}

We also give the asymptotic results for $p_{_{PS}}(t)$ and $p(t)$ for the different scales of $\rho$ and $t$.

\begin{theorem}\label{th7}
The unconditional sojourn time density for the PS model and waiting time density for the ROS model have the following asymptotic expansions:

\begin{enumerate}
\item $\rho$ fixed with $t\to\infty$
\begin{equation}\label{S2_th6_1}
p_{_{PS}}(t)=\frac{p(t)}{1-e^{-\rho}}\sim\frac{\sqrt{\rho+4}-\sqrt{\rho}}{2\sqrt{\rho+4}}\,e^{-1-\rho}\,e^{\rho/\Lambda_0}\,e^{-t}\,e^{\Lambda_0\,t}.
\end{equation}

\item $\rho\to 0$
\begin{enumerate}

\item $t=\zeta/\rho=O(\rho^{-1})$
\begin{equation}\label{S2_th6_21}
p_{_{PS}}(t)\sim \rho^{-1}\,p(t)\sim\frac{1}{2}\,e^{-1}\,e^{-\zeta/2}\,e^{-(1-\sqrt{\rho})t}.
\end{equation}

\item $t=\omega/\sqrt{\rho}=O(\rho^{-1/2})$
\begin{equation}\label{S2_th6_22}
p_{_{PS}}(t)\sim \rho^{-1}\,p(t)\sim e^{-t}\,Q_0(\omega),
\end{equation}
where
$$Q_0(\omega)=\sum_{m=1}^\infty \frac{m^{m-1}}{m!}\,e^{-m}\,\cosh(\omega/\sqrt{m}).$$

\item $t=O(1)$
\begin{equation}\label{S2_th6_23}
p_{_{PS}}(t)=e^{-t}\Big[1+\frac{\rho}{4}\,(t^2-2)+O(\rho^2)\Big].
\end{equation}
\begin{equation}\label{S2_th6_23_p}
p(t)=\rho\,e^{-t}\Big[1+\frac{\rho}{4}\,(t^2-4)+O(\rho^2)\Big].
\end{equation}
\end{enumerate}

\item $\rho\to\infty$ with $t=T\,\rho=O(\rho)$
\begin{equation}\label{S2_th6_31}
p_{_{PS}}(t)\sim p(t)\sim \frac{1}{\rho}\,e^{-T}.
\end{equation}

\end{enumerate}
\end{theorem}

For fixed $\rho$ and large $t$, we removed the condition by using the expansion in the region $t/n>1/\Lambda_0$ (i.e., (\ref{S2_th3_5})) in (\ref{S2_p}), thus obtaining (\ref{S2_th6_1}).

For a small traffic intensity $\rho$, (\ref{S2_th6_21}) on the $t=O(\rho^{-1})$ scale is the limiting case of (\ref{S2_th6_1}), as $\rho\to 0$. For the scale $t=O(\rho^{-1/2})$, we used (\ref{S2_th4_2}) in (\ref{S2_p}). Since $\rho$ is small, the $n=0$ term dominates, which leads to (\ref{S2_th6_22}). When $t=O(1)$, using (\ref{S2_th4_3}) in (\ref{S2_p}) and the fact that $e^{-\rho}\sim 1-\rho$ leads to
$$p_{_{PS}}(t)=(1-\rho)\,\Big[p^{(0)}_0(t)+\rho\,p^{(0)}_1(t)+\rho\,p^{(1)}_0(t)+O(\rho^2)\Big],$$
which yields (\ref{S2_th6_23}). We note that if we let $\omega\to 0$ in (\ref{S2_th6_22}), (\ref{S2_th6_22}) reduces to the leading term in (\ref{S2_th6_23}). This indicates that the $t=O(1)$ scale is a special case of the $t=O(\rho^{-1/2})$ scale, for small $\rho$.

In the case $\rho\to\infty$ with $t=O(\rho)$, we used the leading term in (\ref{S2_th5_1}) in (\ref{S2_p}) and noticed that the infinite sum concentrates near $n=\rho$ (i.e., $N=1$), which led to (\ref{S2_th6_31}). In fact this result is uniform on both the $t=O(\rho)$ and $t=O(\rho^{-1})$ time scales, for large $\rho$.

\section{Derivations of the exact representations}

We first derive the spectral representation (\ref{S2_th1_p}) of the conditional sojourn time density. Consider the equation (\ref{S2_recu2}) and assume that $p_n(t)$ has the form $p_n(t)=e^{\nu\,t}\,\phi(n)$. Then $\phi(n)$ satisfies the recurrence equation
\begin{equation}\label{S3_recu}
(\nu+1)(n+1)\phi(n)=n\,\phi(n-1)-\rho\,\phi(n)+\rho\,\phi(n+1).
\end{equation}
We define the exponential generating function $G(z)$ by
\begin{equation}\label{S3_Gdefine}
G(z)=\sum_{n=0}^\infty\frac{z^n}{n!}\,\phi(n).
\end{equation}
Then by (\ref{S3_recu}), $G(z)$ satisfies
\begin{equation}\label{S3_Gode}
\big[(\nu+1)\,z-\rho\big]\,G'(z)+(\nu+1+\rho-z)\,G(z)=0.
\end{equation}
Here we assumed that $n\,\phi(n-1)$ is finite as $n\to 0$. Solving (\ref{S3_Gode}), we have
\begin{equation}\label{S3_Gresult}
G(z)=C\,\left(1-\frac{\nu+1}{\rho}\,z\right)^{-R_0-1}\,\exp\Big(\frac{z}{\nu+1}\Big),\quad R_0=\frac{\rho\,\nu}{(1+\nu)^2},
\end{equation}
where $C=G(0)=\phi(0)$. Without loss of generality, we let $\phi(0)=C=1$.

To avoid $\phi(n)$ growing like $n!$ as $n\to\infty$, $G(z)$ must be an entire function of $z$, so that $-R_0-1$ must be a non-negative integer. The eigenvalues $\nu$ thus satisfy the quadratic equation $R_0=-m$, $m=1,2,\ldots$, which leads to the two sets of eigenvalues given by (\ref{S2_th1_nu}) and (\ref{S2_th1_nutilde}). We denote by $\phi_m(n,\nu_m)$ and $\phi_m(n,\widetilde{\nu}_m)$ the eigenfunctions corresponding to the eigenvalues $\nu_m$ and $\widetilde{\nu}_m$, respectively. Then for any eigenvalue $\nu_m$, using (\ref{S3_Gresult}) we have
\begin{eqnarray}\label{S3_Gm}
G_m(z)&=&\left(1-\frac{\nu_m+1}{\rho}\,z\right)^{m-1}\,\exp\Big(\frac{z}{\nu_m+1}\Big)\\
&=& \sum_{l=0}^{m-1}\binom{m-1}{l}\left(-\frac{\nu_m+1}{\rho}\right)^{m-1}\,z^{m-1}\,\sum_{k=0}^\infty\frac{z^k}{k!\,(\nu_m+1)^k}.\nonumber
\end{eqnarray}
Thus, from (\ref{S3_Gdefine}) we obtain $\phi_m(n,\nu_m)$ as
\begin{equation*}\label{S3_phimsum}
\phi_m(n,\nu_m)=\sum_{l=0}^{\min(n,m-1)}\binom{m-1}{l}\binom{n}{l}\,l!\,\frac{(\nu_m+1)^{2l-n}}{(-\rho)^l}.
\end{equation*}
We note that as $n\to\infty$, the $l=0$ term dominates and $\phi_m(n,\nu_m)$ is asymptotically given by
\begin{eqnarray}\label{S3_phim_asym}
\phi_m(n,\nu_m)&\sim&\frac{n!}{(n-m+1)!}\,\frac{(\nu_m+1)^{2(m-1)-n}}{(-\rho)^{m-1}}\nonumber\\
&\sim&\frac{n^{m-1}\,(\nu_m+1)^{2m-n-2}}{(-\rho)^{m-1}},\quad n\to\infty.
\end{eqnarray}

Alternately, we can use the generating function of the generalized Laguerre polynomial (see \cite{MA})
\begin{equation}\label{S3_Laguerre}
(1+w)^\alpha\,e^{-w\,x}=\sum_{n=0}^\infty L_n^{(\alpha-n)}(x)\,w^n.
\end{equation}
Comparing (\ref{S3_Gm}) to (\ref{S3_Laguerre}) we see that in our problem, $w=-(\nu_m+1)\,z/\rho$, $x=\rho/(\nu_m+1)^2$ and $\alpha=m-1$. Thus, by the definition of $G$ in (\ref{S3_Gdefine}), we have another representation for $\phi_m(n,\nu_m)$, as
\begin{equation}\label{S3_phimLag}
\phi_m(n,\nu_m)=n!\,\left(-\frac{\nu_m+1}{\rho}\right)^n\,L_n^{(m-1-n)}\left(\frac{\rho}{(\nu_m+1)^2}\right),
\end{equation}
which is (\ref{S2_th1_phi}) with $\nu=\nu_m$. By a similar calculation, we find that the eigenfunctions $\phi_m(n,\widetilde{\nu}_m)$, which correspond to the eigenvalues $\widetilde{\nu}_m$, also satisfy (\ref{S2_th1_phi}) with $\nu=\widetilde{\nu}_m$.

Thus, we can express the conditional sojourn time density as the spectral representation in (\ref{S2_th1_p}), with only the two coefficient sequences $C_m(\nu_m)$ and $C_m(\widetilde{\nu}_m)$ to be determined.

To determine these coefficients, we first obtain an orthogonality relation for the eigenfunctions. Since all of the eigenfunctions satisfy (\ref{S3_recu}), we consider any two eigenfunctions $\phi_m(n,\nu_m)$ and $\phi_{m'}(n,\nu_{m'})$ ($m\ne m'$), which satisfy
\begin{equation}\label{S3_phim}
\big[(\nu_m+1)(n+1)-\rho\big]\,\phi_m(n,\nu_m)=n\,\phi_m(n-1,\nu_m)+\rho\,\phi_m(n+1,\nu_m)
\end{equation}
and
\begin{equation}\label{S3_phim'}
\big[(\nu_{m'}+1)(n+1)-\rho\big]\,\phi_{m'}(n,\nu_{m'})=n\,\phi_{m'}(n-1,\nu_{m'})+\rho\,\phi_{m'}(n+1,\nu_{m'}).
\end{equation}
We multiply (\ref{S3_phim}) by $\rho^n\,\phi_{m'}(n,\nu_{m'})/n!$ and (\ref{S3_phim'}) by $\rho^n\,\phi_{m}(n,\nu_{m})/n!$, subtract one equation from the other, and sum over $n\geq 0$. This leads to
\begin{equation*}
(\nu_m-\nu_{m'})\sum_{n=0}^\infty\frac{n+1}{n!}\,\rho^n\,\phi_m(n,\nu_m)\,\phi_{m'}(n,\nu_{m'})=0.
\end{equation*}
Since $\nu_m\ne \nu_{m'}$, we obtain the orthogonality relation
\begin{equation}\label{S3_orth}
\sum_{n=0}^\infty\frac{n+1}{n!}\,\rho^n\,\phi_m(n,\nu_m)\,\phi_{m'}(n,\nu_{m'})=0.
\end{equation}

By the spectral representation (\ref{S2_th1_p}) and the initial condition (\ref{S2_pn0}), we must have
\begin{equation*}
\sum_{m=1}^\infty C_m(\nu_m)\,\phi_m(n,\nu_m)+\sum_{m=1}^\infty C_m(\widetilde{\nu}_m)\,\phi_m(n,\widetilde{\nu}_m)=\frac{1}{n+1}.
\end{equation*}
Using (\ref{S3_orth}), we can easily show that, for any eigenvalue $\nu_m$ or $\widetilde{\nu}_m$,
\begin{equation}\label{S3_Cm}
C_m(\nu_m)=\frac{\sum_{n=0}^\infty\rho^n\,\phi_m(n,\nu_m)/n!}{\sum_{n=0}^\infty(n+1)\,\rho^n\,\phi_m^2(n,\nu_m)/n!}.
\end{equation}
Using the generating function (\ref{S3_Gm}), the numerator in (\ref{S3_Cm}) is
\begin{equation}\label{S3_nume}
\sum_{n=0}^\infty\frac{\rho^n}{n!}\,\phi_m(n,\nu_m)=G_m(\rho)=(-\nu_m)^{m-1}\,\exp\Big(\frac{\rho}{\nu_m+1}\Big).
\end{equation}

To determine the denominator in (\ref{S3_Cm}), we let $G(z,\mu)$ be a solution of (\ref{S3_recu}) with $\nu=\mu$, whose generating function is given by
\begin{equation*}
G(z,\mu)=\left(1-\frac{\mu+1}{\rho}\,z\right)^{-R-1}\exp\Big(\frac{z}{\mu+1}\Big),\quad R=\frac{\rho\,\mu}{(1+\mu)^2}.
\end{equation*}
Here $\mu$ is not necessarily an eigenvalue, and we assume that $R<-1$ and that $R+1$ is not an integer. Thus, $\phi(n,\mu)$, which corresponds to the non-eigenvalue $\mu$, can be represented as the Cauchy integral
\begin{equation}\label{S3_phimu}
\phi(n,\mu)=\frac{n!}{2\pi i}\oint_\mathcal{C}\,\left(1-\frac{\mu+1}{\rho}\,z\right)^{-R-1}\exp\Big(\frac{z}{\mu+1}\Big)\,z^{-n-1}\,dz,
\end{equation}
where the contour $\mathcal{C}$ is a small circle in the complex $z$-plane centered at the origin. Since $-R-1>0$ and $R$ is not an integer, $z=\rho/(\mu+1)$ is a branch point of the integrand in (\ref{S3_phimu}). We also have the binomial expansion
\begin{equation*}
\left(1-\frac{\mu+1}{\rho}\,z\right)^{-R-1}=\sum_{k=0}^\infty\frac{\Gamma(R+1+k)}{\Gamma(R+1)\,k!}\Big(\frac{\mu+1}{\rho}\Big)^k\,z^k.
\end{equation*}
Then from (\ref{S3_phimu}), we find that as $n\to\infty$, $\phi(n,\mu)$ is asymptotically given by
\begin{eqnarray}\label{S3_phimu_asym}
\phi(n,\mu)&\sim&\frac{\Gamma(R+1+n)}{\Gamma(R+1)}\Big(\frac{\mu+1}{\rho}\Big)^n \exp\Big(\frac{\rho}{(\mu+1)^2}\Big)\nonumber\\
&\sim&\frac{n!\,n^R}{\Gamma(R+1)}\Big(\frac{\mu+1}{\rho}\Big)^n \exp\Big(\frac{\rho}{(\mu+1)^2}\Big),\quad n\to\infty.
\end{eqnarray}

Since $\phi(n,\mu)$ satisfies (\ref{S3_recu}), we have
\begin{equation}\label{S3_phimu2}
\big[(\mu+1)(n+1)-\rho\big]\,\phi(n,\mu)=n\,\phi(n-1,\mu)+\rho\,\phi(n+1,\mu).
\end{equation}
We multiply (\ref{S3_phim}) by $\rho^n\,\phi(n,\mu)/n!$ and (\ref{S3_phimu2}) by $\rho^n\,\phi_{m}(n,\nu_{m})/n!$, subtract one equation from the other, and sum over $0\leq n\leq K$, which yields
\begin{eqnarray}\label{S3_phim+mu}
&&(\nu_m-\mu)\sum_{n=0}^K\frac{n+1}{n!}\,\rho^n\,\phi_m(n,\nu_m)\,\phi(n,\mu)\nonumber\\
&&=\frac{\rho^{K+1}}{K!}\,\big[\phi_m(K+1,\nu_m)\,\phi(K,\mu)-\phi_m(K,\nu_m)\,\phi(K+1,\mu)\big].
\end{eqnarray}
We let $K\to\infty$ and use (\ref{S3_phim_asym}) and (\ref{S3_phimu_asym}), which shows that
$$\phi_m(K+1,\nu_m)\,\phi(K,\mu)=O(K^{m+R-1}),$$
and $$\phi_m(K,\nu_m)\,\phi(K+1,\mu)=O(K^{m+R}).$$
Thus, the second term inside the bracket in the right-hand side of (\ref{S3_phim+mu}) dominates the first and we obtain, after dividing both sides by $\nu_m-\mu$ and expanding for $K\to\infty$,
\begin{eqnarray}\label{S3_sumK}
&&\quad\sum_{n=0}^K\frac{n+1}{n!}\,\rho^n\,\phi_m(n,\nu_m)\,\phi(n,\mu)\nonumber\\
&&\sim\frac{(-1)^{m}\,K^{m+R}\,(\mu+1)^{K+1}\,(\nu_m+1)^{2m-K-2}}{(\nu_m-\mu)\,\Gamma(R+1)\,\rho^{m-1}}\exp\Big(\frac{\rho}{(\mu+1)^2}\Big).
\end{eqnarray}
Next we let $\mu\to\nu_m$ in (\ref{S3_sumK}), so that $R=R(\mu)\to R(\nu_m)=-m$. By using the Laurent expansion of $\Gamma(\cdot)$ near a pole
\begin{equation*}
\Gamma(R+1)\sim\frac{(-1)^{m-1}}{(m-1)!\,(R+m)}\;\;\;\textrm{as   }R\to -m,
\end{equation*}
and then l'H\^{o}pital's rule, we find that
\begin{equation}\label{S3_limmu}
\lim_{\mu\to\nu_m}\Big[(\nu_m-\mu)\,\Gamma(R+1)\Big]=\frac{(-1)^{m}}{(m-1)!\,R'(\nu_m)}=\frac{(-1)^m\,(1+\nu_m)^3}{(m-1)!\,\rho\,(1-\nu_m)}.
\end{equation}
Thus, by using (\ref{S3_sumK}) and (\ref{S3_limmu}) and noting that $\rho\,\nu_m/(\nu_m+1)^2=-m$, we let $K\to\infty$ and obtain
\begin{eqnarray}\label{S3_deno}
&&\sum_{n=0}^\infty\frac{(n+1)}{n!}\,\rho^n\,\phi_m^2(n,\nu_m)\nonumber\\
&&\quad\quad=\frac{m!}{m^{m-1}}\,(-\nu_m)^{m-2}\,(1-\nu_m)\exp\Big(\frac{\rho}{(\nu_m+1)^2}\Big),
\end{eqnarray}
which determines the denominator in (\ref{S3_Cm}). Using (\ref{S3_nume}) and (\ref{S3_deno}) in (\ref{S3_Cm}), we obtain (\ref{S2_th1_C}) with $\nu=\nu_m$. By the same argument, we find that (\ref{S2_th1_C}) is also true for the eigenvalues $\widetilde{\nu}_m$. This completes the derivation of Theorem 2.1.

In the remainder of this section, we use a discrete Green's function to derive (\ref{S2_th2_phat}). Consider the recurrence equation (\ref{S2_recu3}). The discrete Green's function $\mathcal{G}(\theta;n,l)$ satisfies
\begin{eqnarray}\label{S3_dgf}
&&\rho\,\mathcal{G}(\theta;n+1,l)-[(n+1)(1+\theta)+\rho]\,\mathcal{G}(\theta;n,l)\nonumber\\
&&\quad\quad\quad\quad\quad\quad+n\,\mathcal{G}(\theta;n-1,l)=-\delta(n,l),\quad (n,l\geq 0)
\end{eqnarray}
where $\delta(n,l)=1{\{n=l\}}$ is the Kronecker delta. To construct the Green's function we need two linearly independent solutions to
\begin{equation}\label{S3_dgf_Homo}
\rho\,G(\theta;n+1,l)-[(n+1)(1+\theta)+\rho]\,G(\theta;n,l)+n\,G(\theta;n-1,l)=0,
\end{equation}
which is the homogeneous version of (\ref{S3_dgf}).

We seek solutions of (\ref{S3_dgf_Homo}) of the form
$$ G_n=\int_\mathcal{D}z^ng(z)dz,$$
where the function $g(z)$ and the path $\mathcal{D}$ of integration in the complex $z$-plane are to be determined. Using the above form in (\ref{S3_dgf_Homo}) and integrating by parts yields
\begin{eqnarray}\label{S3_IBP}
&&z^ng(z)\big[(1+\theta)\,z-1\big]\Big|_\mathcal{D}\nonumber\\
&&\quad -\int_\mathcal{D}z^n\Big\{\big[(1+\theta)\,z-1\big]g'(z)+\rho\,(z-1)\,g(z)\Big\}dz=0.
\end{eqnarray}
The first term represents contributions from the endpoints of the contour $\mathcal{D}$.

If (\ref{S3_IBP}) is to hold for all $n$ the integrand must vanish, so that $g(z)$ must satisfy the differential equation
\begin{equation}\label{S3_gz}
\big[(1+\theta)\,z-1\big]g'(z)+\rho\,(z-1)\,g(z)=0,
\end{equation}
and thus
$$ g(z)=\Big(z-\frac{1}{1+\theta}\Big)^r\,\exp\Big(-\frac{\rho}{1+\theta}\,z\Big),\;\;r=\frac{\rho\,\theta}{(1+\theta)^2}.$$

If the path of integration $\mathcal{D}$ is chosen as the segment $[0,1/(1+\theta)]$ of the real axis, then (\ref{S3_IBP}) is satisfied for $n\geq 1$. Thus, we obtain $G_n$ as in (\ref{S2_th2_G}). We note that $G_n$ decays as $n\to\infty$, and by scaling $z=(1-y/n)/(1+\theta)$ and using the Laplace method, we find that $G_n$ is asymptotically given by
\begin{equation}\label{S3_Gasymp}
G_n\sim\frac{\Gamma(r+1)}{n^{r+1}\,(1+\theta)^{n+r+1}}\,e^{-r/\theta},\quad n\to\infty.
\end{equation}
However, $G_n$ becomes infinite as $n\to-1$, and $n\,G_{n-1}$ goes to a nonzero limit as $n\to 0$. Thus $G_n$ is not an acceptable solution to (\ref{S3_dgf_Homo}) at $n=0$.

To construct a second solution to (\ref{S3_dgf_Homo}), we consider another path of the integration, the real interval $[1/(1+\theta),\,\infty)$. Then (\ref{S3_IBP}) is again satisfied. Thus, we have another solution of (\ref{S3_dgf_Homo}), $H_n$, which is given by (\ref{S2_th2_H}). $H_n$ is finite as $n\to-1$, but grows as $n\to\infty$. By scaling $z=n\,w=O(n)$ in the integrand of (\ref{S2_th2_H}) we find that $H_n$ grows roughly like $n!$ for $n$ large; more precisely
\begin{equation}\label{S3_Hasymp}
H_n\sim n!\,n^r\,\Big(\frac{1+\theta}{\rho}\Big)^{n+r+1},\quad n\to\infty.
\end{equation}

Thus, the discrete Green's function can be represented by
\begin{eqnarray}\label{S3_G_step}
\mathcal{G}(\theta;n,l) = \left\{ \begin{array}{ll}
H_l\,G_n\,\mathcal{G}_0 & \textrm{if $n\geq l$}\\
G_l\,H_n\,\mathcal{G}_0 & \textrm{if $0\leq n<l$},
\end{array} \right.
\end{eqnarray}
which has acceptable behavior both at $n=0$ and as $n\to\infty$. Here $\mathcal{G}_0$ depends only upon $\theta$ and $l$.

To determine $\mathcal{G}_0$, we let $n=l$ in (\ref{S3_dgf_Homo}) and use the fact that both $G_l$ and $H_l$ satisfy (\ref{S3_dgf_Homo}) with $n=l$. Then we can infer a simple difference equation for the discrete Wronskian $G_l\,H_{l+1}-G_{l+1}\,H_l$, whose solution we write as
\begin{equation}\label{S3_G1}
G_l\,H_{l+1}-G_{l+1}\,H_l=\frac{l!}{\rho\,^l\,\mathcal{G}_1},
\end{equation}
where $\mathcal{G}_1=\mathcal{G}_1(\theta)$ depends upon $\theta$ only. Then using (\ref{S3_G_step}) in (\ref{S3_dgf_Homo}) with $n=l$ shows that $\mathcal{G}_0$ and $\mathcal{G}_1$ are related by $\mathcal{G}_0=\rho^{l-1}\,\mathcal{G}_1/l!$.

Letting $l\to\infty$ in (\ref{S3_G1}) and using the asymptotic results in (\ref{S3_Gasymp}) and (\ref{S3_Hasymp}), we determine $\mathcal{G}_1$ and then obtain $\mathcal{G}_0$ as
$$\mathcal{G}_0=\frac{\rho^{r+l+1}}{l!\,\Gamma(r+1)\,(1+\theta)}\,e^{r/\theta}.$$

Then, we multiply (\ref{S3_dgf}) by the solution $\widehat{p}_l(\theta)$ to (\ref{S2_recu3}) and sum over all $l\geq 0$. After some manipulation this yields
$$\widehat{p}_n(\theta)=\sum_{l=0}^\infty\mathcal{G}(\theta;n,l),$$
which is equivalent to (\ref{S2_th2_phat}). Taking the inverse Laplace transform gives the conditional sojourn time density $p_n(t)$ as the contour integral
\begin{equation}\label{S3_p}
p_n(t)=\frac{1}{2\pi i}\int_{Br}\widehat{p}_n(\theta)e^{\theta t}d\theta,
\end{equation}
where $Br$ is a vertical contour in the complex $\theta$-plane, with $\Re(\theta)\geq 0$. The form in (\ref{S3_p}) is more useful than the spectral representation for obtaining asymptotic results in various limits, such as $n,\,t$ simultaneously large.

To compute the first two moments $\mathcal{M}_n$ and $\mathcal{S}_n$, we could expand $\widehat{p}_n(\theta)$ about $\theta=0$, but it is easier to derive simple difference equations for the moments directly from (\ref{S2_recu2}). By multiplying (\ref{S2_recu2}) by $t$ and integrating from $t=0$ to $t=\infty$ we obtain
\begin{equation*}\label{S3_Mn}
\rho\,\mathcal{M}_{n+1}-(n+1+\rho)\,\mathcal{M}_n+n\,\mathcal{M}_{n-1}=-(n+1).
\end{equation*}
This equation has the linear solution in (\ref{S2_th2_Mn}). Similarly, for the second moment we obtain
\begin{equation}\label{S3_Sn}
\rho\,\mathcal{S}_{n+1}-(n+1+\rho)\,\mathcal{S}_n+n\,\mathcal{S}_{n-1}=-2(n+1)\,\mathcal{M}_n,
\end{equation}
whose solution is given by (\ref{S2_th2_Sn}).

\section{Asymptotic results for fixed $\rho$ and $\rho\to 0$}

We first assume that the traffic intensity $\rho$ is fixed. We sketch the main points in deriving Theorem 2.3. We first consider $n,t\to\infty$ with $n>t$ and use the result in (\ref{S2_th2_phat}). To obtain a two term approximation, we need the correction terms in the approximations in (\ref{S3_Gasymp}) and (\ref{S3_Hasymp}), which are given by
\begin{equation}\label{S4_Gasymp}
G_n=\frac{e^{-r/\theta}}{n^{r+1}\,(1+\theta)^{n+r+1}}\,\Big[\Gamma(r+1)+\frac{1}{n}\Big(\frac{r}{\theta}\,\Gamma(r+2)-\frac{1}{2}\,\Gamma(r+3)\Big)+O(n^{-2})\Big]
\end{equation}
and
\begin{equation}\label{S4_Hasymp}
H_n= n!\,n^r\,\Big(\frac{1+\theta}{\rho}\Big)^{n+r+1}\left[1-\frac{r^2}{n\,\theta}+O(n^{-2})\right].
\end{equation}
From (\ref{S4_Gasymp}) and (\ref{S4_Hasymp}), we note that the first term in (\ref{S2_th2_phat}) dominates the second, and thus the Laplace transform is asymptotically given by
\begin{eqnarray}
\widehat{p}_n(t)&\sim& M\,G_n\sum_{l=0}^n\,\frac{\rho^l}{l!}\,H_l\label{S4_1_phat}\\
&\sim&\frac{1+\theta}{\rho\,\theta}\sum_{l=0}^n (1+\theta)^{l-n}\Big(\frac{l}{n}\Big)^r\,\frac{1}{n}\,\Big[r+\frac{A}{n}-\frac{r^2}{l\,\theta}\Big]\label{S4_1_sum}\\
&\sim&\Big[\frac{1}{1+\theta}+\frac{(1+\theta)\,A}{\rho\,\theta\,n}\Big]\int_0^1(1+\theta)^{-n\,y}(1-y)^r\,dy\nonumber\\
&&-\frac{\rho}{(1+\theta)^3\,n}\int_0^1(1+\theta)^{-n\,y}(1-y)^{r-1}\,dy\label{S4_1_int}\\
&\equiv& I_1+I_2.\nonumber
\end{eqnarray}
Here
$$A=A(\theta)\equiv\frac{r\,(r+1)}{2\theta}\big[2r-\theta\,(r+2)\big]$$
and we used the Euler-Maclaurin summation formula to approximate the sums in (\ref{S4_1_sum}) by integrals. By scaling $\theta=s/n=O(1/n)$ and noting that
$$ \frac{(1+\theta)\,A}{\rho\,\theta}=\rho-1+\frac{s}{2n}(2\rho^2-9\rho+2)+O(n^{-2})$$
and $r\sim \rho\,s/n$, the first term, $I_1$, in (\ref{S4_1_int}) becomes
\begin{eqnarray*}
I_1&\sim& \Big(1+\frac{\rho-1}{n}\Big)\int_0^1 e^{-s\,y}\,dy\\
&=& \Big(1+\frac{\rho-1}{n}\Big)\frac{1-e^{-s}}{s}.
\end{eqnarray*}
Thus, taking the inverse Laplace transforms of $I_1$ and $I_2$ yields
\begin{eqnarray}\label{S4_invI}
\mathcal{L}^{-1}(I_1)&\sim& \Big(1+\frac{\rho-1}{n}\Big)\,\frac{1}{n}\,\textrm{Heaviside}\Big(1-\frac{t}{n}\Big)\nonumber\\
&=&\frac{1}{n}+\frac{\rho-1}{n^2}\quad\quad (n>t)
\end{eqnarray}
and
\begin{eqnarray}\label{S4_invII}
\mathcal{L}^{-1}(I_2)&\sim&-\frac{\rho}{n^2}\,\frac{1}{2\pi i}\int_{Br}e^{s\,t/n}\left[\int_0^1\frac{e^{-s\,y}}{(1-y)^{1-\rho s/n}}\,dy\right]\,ds\nonumber\\
&\sim&-\frac{\rho}{n^2}\int_0^1\frac{\delta(t/n-y)}{1-y}\,dy\nonumber\\
&=&-\frac{\rho}{n\,(n-t)}\quad\quad (n>t).
\end{eqnarray}
We note that in deriving (\ref{S4_invII}), we changed the order of integration and used $d\theta=n^{-1}ds$. Then (\ref{S4_invI}) and (\ref{S4_invII}) lead to (\ref{S2_th3_1}).

This analysis suggests that $p_n(t)$ is approximately zero in the range $t/n>1$. We shall show that in this sector the density is exponentially small. Before doing this, we first investigate the transition region, where $t\approx n$.

Thus, we consider $n,t\to\infty$ with $n/t=1+\Delta\,t^{-1/2}=1+O(t^{-1/2})$. We can still use (\ref{S4_1_phat}) but now scale $l=y\,\sqrt{n}=O(\sqrt{n})$, and approximate the sum by
\begin{equation}\label{S4_2_phat}
\widehat{p}_n(\theta)\sim\frac{1}{(1+\theta)^{n+1}\,(\sqrt{n}\,)^{r+1}}\int_0^\infty(1+\theta)^{\sqrt{n}\,y}\,y^r\,dy.
\end{equation}
Scaling $\theta=\varpi/\sqrt{n}=O(1/\sqrt{n})$, and noting that
\begin{equation*}
(1+\theta)^{-n-1}\sim e^{\varpi^2/2-\sqrt{n}\,\varpi}
\end{equation*}
and $(1+\theta)^{\sqrt{n}\,y}\sim e^{\varpi\,y}$, the inverse Laplace transform leads to
\begin{eqnarray}\label{S4_2_pnt}
p_n(t)&\sim&\frac{1}{n}\,\frac{1}{2\pi i}\int_{Br}e^{\varpi\,t/\sqrt{n}}\,e^{\varpi^2/2-\sqrt{n}\,\varpi}\,\int_0^\infty e^{\varpi\,y}\,dy\,d\varpi\nonumber\\
&=&\frac{1}{n}\int_0^\infty\frac{1}{2\pi i}\int_{Br}e^{\varpi^2/2}\,e^{(-\sqrt{n}+y+t/\sqrt{n})\,\varpi}d\varpi\,dy.
\end{eqnarray}
Note that in this range of $(n,t)$,
\begin{equation*}
\frac{t}{\sqrt{n}}-\sqrt{n}=\Big(1-\frac{n}{t}\Big)\,\frac{t}{\sqrt{n}}=-\Delta\frac{\sqrt{t}}{\sqrt{n}}\sim-\Delta=O(1).
\end{equation*}
Then by using the identity
\begin{equation*}
\frac{1}{2\pi i}\int_{Br}e^{C_0\varpi^2+C_1\varpi}d\varpi=\frac{1}{2\sqrt{\pi\,C_0}}\exp\Big(-\frac{C_1^2}{4C_0}\Big)
\end{equation*}
and noting that $\Delta=(n-t)/\sqrt{t}$, we explicitly evaluate the integral over $\varpi$ in (\ref{S4_2_pnt}) and obtain (\ref{S2_th3_2}).

Now we consider $n,\;t\to\infty$ with $t>n$. We rewrite (\ref{S2_th2_phat}) as
\begin{equation}\label{S4_3_phat}
\widehat{p}_n(\theta)=M\,G_n\sum_{l=0}^\infty\frac{\rho^l}{l!}\,H_l\;+\;M\sum_{l=n+1}^\infty\frac{\rho^l}{l!}\,(H_n\,G_l-H_l\,G_n).
\end{equation}
The first sum can be calculated exactly by using (\ref{S2_th2_H}) and making the transformation $z=1/(1+\theta)+y$, which yields
\begin{equation}\label{S4_3_sumH}
\sum_{l=0}^\infty \frac{\rho^l}{l!}\,H_l=e^r\,\Gamma(r+1)\left(-\frac{1+\theta}{\rho\,\theta}\right)^{r+1}.
\end{equation}
The result in (\ref{S4_3_sumH}) holds for $\theta<0$ and $\theta>\theta_p=-1+[-\rho+\sqrt{\rho^2+4\rho}]/2$, since $\Gamma(r+1)$ has a simple pole at $\theta=\theta_p$. We note that $\theta_p=\nu_1$, which is the first eigenvalue in (\ref{S2_th1_nu}). The second sum in (\ref{S4_3_phat}) is negligible in view of (\ref{S3_Gasymp}), (\ref{S3_Hasymp}) and (\ref{S4_3_sumH}), and the fact $\theta<0$ on this scale. Using (\ref{S2_th2_M}), (\ref{S3_Gasymp}) and (\ref{S4_3_sumH}) in the first sum of (\ref{S4_3_phat}), then taking the inverse Laplace transform, we have
\begin{equation}\label{S4_3_pnt}
p_n(t)\sim\frac{1}{2\pi i}\int_{Br}h(\theta)\,e^{t\,f(\theta)}\,d\theta,
\end{equation}
where $f(\theta)=\theta-\log(1+\theta)\,n/t$ and
\begin{equation*}
h(\theta)=\frac{\Gamma(r+1)\,e^r}{(1+\theta)\,(-\theta)^{r+1}\,n^{r+1}}.
\end{equation*}
For $t\to\infty$ and $n/t$ fixed we evaluate (\ref{S4_3_pnt}) by the saddle point method.
There is a saddle point at $\theta=\theta_s\equiv n/t-1<0$, which satisfies $f'(\theta)=0$. Hence, using the saddle point method gives
\begin{equation*}
p_n(t)\sim\frac{h(\theta_s)}{\sqrt{2\pi\,t\,f''(\theta_s)}}\,e^{t\,f(\theta_s)}
\end{equation*}
and this leads to (\ref{S2_th3_3}), where
\begin{equation*}
r_\ast=r(\theta_s)=\frac{\rho\,\theta_s}{(1+\theta_s)^2}=\rho\,\frac{t^2}{n^2}\,\Big(\frac{n}{t}-1\Big).
\end{equation*}
This analysis indicates that (\ref{S2_th3_3}) only holds for $n,\;t\to\infty$ with $n/t<1$ and $\theta_s>\theta_p$, so that $n/t=1+\theta_s>1+\theta_p=\Lambda_0$.

There is a transition region where $n/t=\Lambda_0+\Lambda/\sqrt{t}$ with $\Lambda=O(1)$. We still use (\ref{S4_3_pnt}) and note that $h(\theta)$ has a simple pole at $\theta=\theta_p$, and the saddle point $\theta_s$ of $f(\theta)$ is now close to $\theta_p$. We expand the integrand in (\ref{S4_3_pnt}) about $\theta=\theta_p$ using
\begin{equation*}
h(\theta)\sim\frac{\sqrt{\rho+4}-\sqrt{\rho}}{2\sqrt{\rho+4}}\,e^{-1}\,\frac{1}{\theta-\theta_p}
\end{equation*}
and
\begin{eqnarray*}
f(\theta)&=&f(\theta_p)+f'(\theta_p)(\theta-\theta_p)+\frac{1}{2}f''(\theta_p)(\theta-\theta_p)^2+\cdots\\
&\sim& (\Lambda_0-1)-\frac{n}{t}\,\log(\Lambda_0)-\frac{\Lambda}{\Lambda_0}\,\frac{\theta-\theta_p}{\sqrt{t}}+\frac{1}{2\Lambda_0}(\theta-\theta_p)^2+\cdots.
\end{eqnarray*}
Here we also used $1+\theta_p=\Lambda_0$ and $n/t\sim\Lambda_0$. By scaling $\theta-\theta_p=S/\sqrt{t}$, (\ref{S4_3_pnt}) asymptotically becomes
\begin{equation*}
p_n(t)\sim\Lambda_0^{-n}\,e^{-t+\Lambda_0\,t}\,\frac{1}{2\pi i}\int_{Br}\frac{1}{S}\,\exp\Big[-\frac{\Lambda}{\Lambda_0}\,S+\frac{1}{2\Lambda_0}\,S^2\Big]dS,
\end{equation*}
where $\Re(S)>0$ on the contour $Br$. Then we use the identity
\begin{equation*}
\frac{1}{2\pi i}\int_{Br}\frac{1}{S}\,e^{-A\,S+B\,S^2/2}dS=\frac{1}{\sqrt{2\pi}}\int_{A/\sqrt{B}}^\infty e^{-u^2/2}du,
\end{equation*}
with $A=\Lambda/\Lambda_0$ and $B=\Lambda_0^{-1}$, to eventually obtain (\ref{S2_th3_4}).

Finally, for the scale $n,\;t\to\infty$ with $n/t<\Lambda_0$, the pole at $\theta=\theta_p$ dominates the asymptotic behavior of $p_n(t)$, and (\ref{S2_th3_5}) is obtained by evaluating the residue at the dominant pole in (\ref{S4_3_pnt}). This concludes the derivation of Theorem 2.3.

If we consider a small traffic intensity, $\rho\to 0$, and scale $t=O(\rho^{-1})$, we obtain seven different asymptotic expressions based on different scalings of the space variable $n$, which are given in case 1 of Theorem 2.4. All seven results can be obtained as limiting cases of the results in Theorem 2.3, by letting $\rho\to 0$ and scaling $n$ appropriately. We omit the derivations here. Note that as $\rho\to 0$ the transition line $n/t=\Lambda_0$ becomes close to the $t$-axis.

For the time scale $t=O(\rho^{-1/2})$, from the spectral representation we note that as $\rho\to 0$ the eigenvalues are
\begin{equation*}
\nu_m=-1+\frac{\sqrt{\rho}}{\sqrt{m}}+O(\rho),\quad \widetilde{\nu}_m=-1-\frac{\sqrt{\rho}}{\sqrt{m}}+O(\rho).
\end{equation*}
Then the eigenfunctions are asymptotically given by
\begin{equation*}
\phi_m(n,\nu_m)\sim(-1)^n\,n!\,m^{-n/2}\,\rho^{-n/2}\,L_n^{(m-1-n)}(m)
\end{equation*}
and $\phi_m(n,\widetilde{\nu}_m)\sim (-1)^n\,\phi_m(n,\nu_m)$.
Thus, all the eigenvalues contribute to $p_n(t)$ on the scale $t=O(\rho^{-1/2})$ and $n=O(1)$, and we obtain (\ref{S2_th4_2}).

Now we consider the scale $n,\;t=O(1)$ with $\rho\to 0$ and use the result in (\ref{S2_th2_phat}). Since $r=O(\rho)$, $G_n$ in (\ref{S2_th2_G}) becomes
\begin{equation}\label{S4_th4_3_G}
G_n\sim\int_0^{\frac{1}{1+\theta}}z^n\,dz=\frac{1}{(n+1)\,(1+\theta)^{n+1}}.
\end{equation}
By scaling $z=\frac{s}{(1+\theta)\,\rho}$ in (\ref{S2_th2_H}), $H_n$ is asymptotically given by
\begin{eqnarray}\label{S4_th4_3_H}
H_n&\sim&\frac{1}{\rho^{n+1}\,(1+\theta)^{n+1}}\int_0^\infty s^n\,\exp\Big(-\frac{s}{(1+\theta)^2}\Big)\,ds\nonumber\\
&=& n!\,\frac{(1+\theta)^{n+1}}{\rho^{n+1}}.
\end{eqnarray}
We also have $M\sim\rho/(1+\theta)$. Using (\ref{S4_th4_3_G}) and (\ref{S4_th4_3_H}) in (\ref{S2_th2_phat}), we find that the first sum dominates the second and we obtain
\begin{eqnarray*}
\widehat{p}_n(\theta)&\sim& M\,G_n\sum_{l=0}^n\frac{\rho^l}{l!}\,H_l\\
&\sim& \frac{1}{(n+1)\,(1+\theta)^{n+2}}\sum_{l=0}^{n}(1+\theta)^{l+1}\\
&=&\frac{1}{n+1}\,\frac{1}{\theta}\,\left[1-\frac{1}{(1+\theta)^{n+1}}\right].
\end{eqnarray*}
Then we invert the Laplace transform over time, letting $\theta=w-1$, which gives
\begin{eqnarray}\label{S4_th4_3_p0}
p_n(t)&\sim&\frac{e^{-t}}{n+1}\,\frac{1}{2\pi i}\int_{Br}\frac{1-w^{-n-1}}{w-1}\,e^{w\,t}\,dw\nonumber\\
&=&\frac{e^{-t}}{n+1}\sum_{l=0}^n\frac{t^l}{l!}.
\end{eqnarray}
This is the leading term $p^{(0)}_n(t)$ in (\ref{S2_th4_3}).

To obtain the second term $p^{(1)}_n(t)$, we need the correction terms of the asymptotic expansions in (\ref{S4_th4_3_G}) and (\ref{S4_th4_3_H}). It is much easier, however, to use a perturbation method to obtain $p^{(1)}_n(t)$. Assume that the conditional sojourn time density has an expansion in powers of $\rho$, as in (\ref{S2_th4_3}). By using the recurrence equation (\ref{S2_recu2}), the leading term $p^{(0)}_n(t)$ satisfies
$$\frac{d\,p^{(0)}_n(t)}{dt}=\frac{n}{n+1}\,p^{(0)}_{n-1}(t)-p^{(0)}_n(t)$$
with the initial condition $p^{(0)}_n(t)=1/(n+1)$. This can be easily solved to regain (\ref{S4_th4_3_p0}). The second term $p^{(1)}_n(t)$ satisfies
\begin{equation*}
\frac{d\,p^{(1)}_n(t)}{dt}=\frac{n}{n+1}\,p^{(1)}_{n-1}(t)-p^{(1)}_n(t)+\frac{1}{n+1}\,p^{(0)}_{n+1}(t)-\frac{1}{n+1}\,p^{(0)}_n(t)
\end{equation*}
with the initial condition $p^{(1)}_n(0)=0$. We set
\begin{equation}\label{S4_th4_3_p1}
p^{(1)}_n(t)=\frac{e^{-t}}{n+1}\,\mathcal{P}_n(t)
\end{equation}
and take the Laplace transform of $\mathcal{P}_n(t)$ over time, with $\widehat{\mathcal{P}}_n(s)=\int_0^\infty\mathcal{P}_n(t)\,e^{-st}dt$. After some simplification, we find that $\widehat{\mathcal{P}}_n(s)$ satisfies the following difference equation:
\begin{equation*}
s\,\widehat{\mathcal{P}}_n(s)-\widehat{\mathcal{P}}_{n-1}(s)=-\frac{(n+1)\,s^{-n-2}-(n+2)\,s^{-n-1}+1}{(n+1)(n+2)(s-1)},\quad n\geq 1
\end{equation*}
with
$$\widehat{\mathcal{P}}_0(s)=\frac{1-s}{2s^3}.$$
After some calculation, we obtain $\widehat{\mathcal{P}}_n(s)$ as
\begin{eqnarray*}
\widehat{\mathcal{P}}_n(s)&=&-\frac{s^{-n-1}}{s-1}+\frac{1}{(n+2)s(s-1)}+\frac{(n+1)s^{-n}}{(n+2)s^2(s-1)}\\
&&+\frac{1}{s^{n+3}}\sum_{l=0}^n\frac{1}{l+2}-\sum_{l=0}^{n-1}\frac{s^{l-n-1}}{l+2}.
\end{eqnarray*}
Then by taking the inverse Laplace transform, and using the relation (\ref{S4_th4_3_p1}), we obtain $p^{(1)}_n(t)$ in (\ref{S2_th4_3}).

\section{Asymptotic results for $\rho\to\infty$}

We shall use a singular perturbation approach to derive the asymptotic approximations for large traffic intensities, $\rho\to\infty$. This method should be useful for analyzing models with more general balking probabilities. We shall also sketch how the asymptotic results in Theorem 2.5 follow from the exact representations in Theorems 2.1 and 2.2.

We first consider the scale $t=T\rho=O(\rho)$ and $n=N\rho=O(\rho)$, and expand $p_n(t)$ in powers of $\rho^{-1}$, as in (\ref{S2_th5_1}). Using (\ref{S2_th5_1}) in the recurrence equation (\ref{S2_recu2}), the leading term $P_0(N,T)$ satisfies
\begin{equation}\label{S5_1_p0_pde}
\frac{\partial P_0}{\partial T}+\frac{N-1}{N}\frac{\partial P_0}{\partial N}=-\frac{1}{N}\,P_0
\end{equation}
with the initial condition $P_0(N,T)=1/N$. We solve this first order PDE by the method of characteristics. The family of characteristics is given by
$$T=N+\log|1-N|+\textrm{constant},$$
where the constant indexes the family. The characteristic $T=N+\log(1-N)$ goes through the origin $(N,T)=(0,0)$, along the parabola $T=-N^2/2$. The general solution to (\ref{S5_1_p0_pde}) is
\begin{equation}\label{S5_1_p0_F}
P_0(N,T)=\frac{1}{N-1}\mathcal{F}\big((N-1)\,e^{N-T}\big).
\end{equation}
Using the initial condition in (\ref{S5_1_p0_F}), we determine the function $\mathcal{F}(\cdot)$ from
$$\mathcal{F}\big((N-1)\,e^{N}\big)=\frac{N-1}{N}.$$
If we denote by $N_\ast=N_\ast(N,T)$ the solution to
\begin{equation}\label{S5_1_nstar}
(N_\ast-1)\,e^{N_\ast}=(N-1)\,e^{N-T},
\end{equation}
$P_0$ in (\ref{S5_1_p0_F}) becomes
\begin{equation}\label{S5_1_p0_nstar}
P_0(N,T)=\frac{1}{N-1}\,\frac{N_\ast-1}{N_\ast}.
\end{equation}
Setting $N_\ast=N-U$ in (\ref{S5_1_nstar}) leads to (\ref{S2_th5_1U}). Thus, (\ref{S5_1_p0_nstar}) can be rewritten as (\ref{S2_th5_1p0}).

Alternately, we can rewrite (\ref{S2_th5_1p0}) more explicitly, in terms of an infinite series. From (\ref{S2_th5_1U}), we let $U=N-1+U_0$, where $U_0=U_0(N,T)=(1-N)\,e^{U-T}$. Then $U_0$ can be expressed in terms of the Lambert W-function (see \cite{Lambert}), which satisfies
$$e^{-U_0}\,U_0=(1-N)\,e^{N-T-1}\equiv z.$$
We use the series expansion of the Lambert W-function to obtain $U_0$ as
$$U_0=\sum_{m=1}^\infty\frac{(-m)^{m-1}}{m!}\,z^m,$$
where the series converges for $|z|<e^{-1}$. Thus, $U$ has the following series expansion
\begin{equation}\label{S5_1_U}
U(N,T)=N-1+\sum_{m=1}^\infty\frac{m^{m-1}}{m!}\,(1-N)^m\,e^{m(N-T-1)},
\end{equation}
which converges for $|1-N|\,e^{N-T}<1$. The series is always convergent for $N\leq 1$, but diverges for $N>1$, if $T<N+\log(N-1)$. For example, if $T=0$ the series converges only for $N<N_c\doteq 1.2784$, where $(N_c-1)\,e^{N_c}=1$.
Using (\ref{S5_1_U}) in (\ref{S2_th5_1p0}), we have an alternate series expression for $P_0(N,T)$:
\begin{equation}\label{S5_1_p0_series}
P_0(N,T)=\frac{\displaystyle\sum_{m=1}^\infty m^{m-1}\,(1-N)^{m-1}\,e^{m(N-T-1)}/m!}{1-\displaystyle\sum_{m=1}^\infty m^{m-1}\,(1-N)^m\,e^{m(N-T-1)}/m!}.
\end{equation}

Now we sketch how to compute the correction term $P_1(N,T)$, which satisfies the following PDE
\begin{equation}\label{S5_1_p1_pde}
\frac{\partial P_1}{\partial T}=\frac{1-N}{N}\,\frac{\partial P_1}{\partial N}-\frac{1}{N}\,P_1+\frac{N+1}{2N}\,\frac{\partial^2 P_0}{\partial N^2}+\frac{N-1}{N^2}\,\frac{\partial P_0}{\partial N}+\frac{1}{N^2}\,P_0,
\end{equation}
with the initial condition $P_1(N,0)=-1/N^2$. This follows from expanding $p_n(0)=1/(n+1)=1/(N\rho+1)$ in powers of $\rho^{-1}$.  We make the substitution
\begin{equation}\label{S5_1_p1_sub}
P_1(N,T)=-\frac{1}{N}\,P_0(N,T)+P_1^\ast(N,T),
\end{equation}
so that $P_1^\ast(N,0)=0$.

We introduce an operator $\mathcal{D}$, which is defined by
\begin{equation}\label{S5_1_L}
\mathcal{D}=\frac{\partial}{\partial T}+\frac{N-1}{N}\,\frac{\partial}{\partial N}+\frac{1}{N}.
\end{equation}
Then by (\ref{S5_1_p0_pde}) we have $\mathcal{D}P_0=0$ and
\begin{equation}\label{S5_1_LpoN}
\mathcal{D}\Big(-\frac{P_0}{N}\Big)=\frac{N-1}{N^3}\,P_0.
\end{equation}
Applying the operator $\mathcal{D}$ to (\ref{S5_1_p1_sub}), and using (\ref{S5_1_p1_pde}) and (\ref{S5_1_LpoN}), we obtain
\begin{equation}\label{S5_1_Lpstar}
\mathcal{D}P_1^\ast=\frac{N+1}{2N}\,\frac{\partial^2 P_0}{\partial N^2}+\frac{N-1}{N^2}\,\frac{\partial P_0}{\partial N}+\frac{1}{N^3}\,P_0.
\end{equation}

We change variables in (\ref{S5_1_Lpstar}) from $(N,T)$ to $(\xi,\eta)$, where $N=\xi+\eta$ and $U=\eta$, with $U=U(N,T)$ given by (\ref{S2_th5_1U}). Thus, by the chain rule, we have
\begin{equation*}
\frac{\partial}{\partial N}=\Big(1-\frac{\partial U}{\partial N}\Big)\,\frac{\partial}{\partial \xi}+\frac{\partial U}{\partial N}\,\frac{\partial}{\partial \eta}
\end{equation*}
and
\begin{equation*}
\frac{\partial}{\partial T}=-\frac{\partial U}{\partial T}\,\frac{\partial}{\partial \xi}+\frac{\partial U}{\partial T}\,\frac{\partial}{\partial \eta}.
\end{equation*}
By implicitly differentiating (\ref{S2_th5_1U}) with respect to $N$ and $T$ we obtain
$$\frac{\partial U}{\partial N}=\frac{\eta}{\xi\,(\xi+\eta-1)},\quad \frac{\partial U}{\partial T}=\frac{\xi-1}{\xi}.$$
Then by (\ref{S5_1_L}), after changing variables, the operator $\mathcal{D}$ can be rewritten as
\begin{equation}\label{S5_1_Lchange}
\mathcal{D}=\frac{\xi+\eta-1}{\xi+\eta}\,\frac{\partial}{\partial \eta}+\frac{1}{\xi+\eta}.
\end{equation}
If we denote $P_1^\ast(N,T)$ by $F(\xi,\eta)$, then by the chain rule we have
\begin{equation}\label{S5_1_F1}
\frac{\partial P_0}{\partial N}=\Big(1-\frac{\partial U}{\partial N}\Big)\,\frac{\partial F}{\partial \xi}+\frac{\partial U}{\partial N}\,\frac{\partial F}{\partial\eta}\equiv F_1(\xi,\eta)
\end{equation}
and
\begin{equation}\label{S5_1_F2}
\frac{\partial^2 P_0}{\partial N^2}=\Big(1-\frac{\partial U}{\partial N}\Big)\,\frac{\partial F_1}{\partial \xi}+\frac{\partial U}{\partial N}\,\frac{\partial F_1}{\partial\eta}\equiv F_2(\xi,\eta).
\end{equation}
Thus, by using (\ref{S5_1_Lpstar})-(\ref{S5_1_F2}) and noting that $N=\xi+\eta$, we have
\begin{eqnarray*}
\mathcal{D}F&=&\frac{\xi+\eta-1}{\xi+\eta}\,\frac{\partial F}{\partial\eta}+\frac{1}{\xi+\eta}\,F\\
&=&\frac{\xi+\eta+1}{2(\xi+\eta)}\,F_2+\frac{\xi+\eta-1}{(\xi+\eta)^2}\,F_1+\frac{1}{(\xi+\eta)^3}\,F.
\end{eqnarray*}
Multiplying the above by $\xi+\eta\, (=N)$ and using (\ref{S5_1_Lchange}) we obtain
\begin{equation}\label{S5_1_F_pde}
\frac{\partial}{\partial\eta}\Big[(\xi+\eta-1)\,F\Big]=\frac{\xi+\eta+1}{2}\,F_2+\frac{\xi+\eta-1}{\xi+\eta}\,F_1+\frac{1}{(\xi+\eta)^2}\,F.
\end{equation}
Solving (\ref{S5_1_F_pde}) with the help of the symbolic computation program MAPLE, then imposing the initial condition $P_1^\ast(N,0)=F(\xi,0)=0$ and using (\ref{S5_1_p1_sub}), we obtain the correction term $P_1(N,T)$ in (\ref{S2_th5_1}) as follows
\begin{eqnarray}\label{S5_1_P1}
P_1(N,T)&=&F(\xi,\eta)\nonumber\\
&=&-\frac{(\xi-1)(2\xi-3)}{2\xi^5}+\frac{3(\xi-1)(2\xi^2-2\xi-3)}{2\xi^5(\xi+\eta-1)}\nonumber\\
&&-\frac{(\xi-1)^2(2\xi^3+2\xi^2-5\xi-15)}{2\xi^5(\xi+\eta-1)^2}-\frac{(\xi-1)^3(2\xi^2+4\xi+3)}{2\xi^5(\xi+\eta-1)^3}\nonumber\\
&&-\frac{2(\xi-1)(2\xi-3)}{\xi^5(\xi+\eta-1)}\log\bigg|\frac{\xi+\eta-1}{\xi-1}\bigg|,
\end{eqnarray}
where $\xi=N-U,\;\eta=U.$

Now we consider some special cases. If $N=1$, then $U\to 0$ by (\ref{S2_th5_1U}). Thus, $\xi\sim 1$, $\eta\sim 0$ and
$$\frac{\xi-1}{\xi+\eta-1}\sim e^{-T}.$$
Then (\ref{S5_1_P1}) reduces to the explicit result
\begin{equation*}
P_1(1,T)=\Big(2T-\frac{9}{2}\Big)\,e^{-T}+8\,e^{-2T}-\frac{9}{2}\,e^{-3T}.
\end{equation*}
We already showed that $P_0(1,T)=e^{-T}$. We relate the explicit result along $N=1$ ($n=\rho$) to the spectral expansion (\ref{S2_th1_p}). From (\ref{S2_th1_p}) we can easily show that only the eigenvalues $\nu_1$, $\nu_2$ and $\nu_3$ are $O(\rho^{-1})$ or $O(\rho^{-2})$ in this limit (the others are $o(\rho^{-2})$). Expanding (\ref{S2_th1_p})-(\ref{S2_th1_phi}) for $\rho\to\infty$ and $N=1$, we obtain
\begin{equation*}
C_1(\nu_1)\,\phi_1(\rho,\nu_1)\,e^{\nu_1\,t}\sim \rho^{-1}\,e^{-T}+\rho^{-2}\,\Big(2T-\frac{9}{2}\Big)\,e^{-T},
\end{equation*}
\begin{equation*}
C_2(\nu_2)\,\phi_2(\rho,\nu_2)\,e^{\nu_2\,t}\sim 8\,\rho^{-2}\,e^{-2T},
\end{equation*}
and
\begin{equation*}
C_3(\nu_3)\,\phi_3(\rho,\nu_3)\,e^{\nu_3\,t}\sim -\frac{9}{2}\,\rho^{-2}\,e^{-3T}.
\end{equation*}
Thus (\ref{S2_th1_p}) agrees precisely with $\rho^{-1}\,P_0(1,T)+\rho^{-2}\,P_1(1,T)$.

If $N=0$ and $T\to\infty$, then $U\to -1$, $\eta\sim -1$, and $\xi\sim 1-e^{-1}\,e^{-T}$. Hence, from (\ref{S5_1_P1}), $P_0(0,T)\sim e^{-1}\,e^{-T}$ and $P_1(0,T)\sim (2T-3)\,e^{-1}\,e^{-T}$.
From (\ref{S2_th1_p}) we obtain for $n=0$, $t=T\rho\to\infty$
\begin{equation*}
C_1(\nu_1)\,\phi_1(0,\nu_1)\,e^{\nu_1\,t}\sim \rho^{-1}\,e^{-1}\,e^{-T}+\rho^{-2}\,(2T-3)\,e^{-1}\,e^{-T}.
\end{equation*}
Again, this agrees with (\ref{S5_1_P1}) and shows that for $N=0$ and $T\gg 1$ ($t\gg \rho$) only the first eigenvalue $\nu_1$ contributes to the expansion of $p_n(t)$.

Next, we consider short time scales, with $n=O(1)$ and $t=\tau/\rho=O(\rho^{-1})$. Expanding the conditional sojourn time density $p_n(t)$ in the form $p_n(t)=\mathcal{Q}_n(\tau)+O(\rho^{-1})$ and using equation (\ref{S2_recu2}), we have
\begin{equation*}
\mathcal{Q}'_n(\tau)=\frac{1}{n+1}\Big[\mathcal{Q}_{n+1}(\tau)-\mathcal{Q}_n(\tau)\Big],
\end{equation*}
with the initial condition $\mathcal{Q}_n(0)=1/(n+1)$. Taking the Laplace transform over the time variable $\tau$ with $\widehat{\mathcal{Q}}_n(s)=\int_0^\infty \mathcal{Q}_n(\tau)\,e^{-\tau s}\,d\tau$, we obtain the following difference equation for $\widehat{\mathcal{Q}}_n(s)$:
\begin{equation}\label{S5_2_Qhat_pde}
\widehat{\mathcal{Q}}_{n+1}(s)-\big[(n+1)\,s+1\big]\,\widehat{\mathcal{Q}}_n(s)=-1.
\end{equation}
Solving (\ref{S5_2_Qhat_pde}) yields
\begin{equation}\label{S5_2_Qhat}
\widehat{\mathcal{Q}}_n(s)=\sum_{j=0}^\infty s^{-j-1}\,\frac{\Gamma(n+1+1/s)}{\Gamma(n+j+2+1/s)}.
\end{equation}
By the inverse Laplace transform, we have
\begin{equation}\label{S5_2_Q_inv}
\mathcal{Q}_n(\tau)=\frac{1}{2\pi i}\int_{Br}\frac{e^{\tau s}}{s}\sum_{j=0}^\infty s^{-j}\,\frac{\Gamma(n+1+1/s)}{\Gamma(n+j+2+1/s)}\,ds.
\end{equation}
Using the identity
$$\int_0^1 t^{x-1}\,(1-t)^{y-1}\,dt=\frac{\Gamma(x)\,\Gamma(y)}{\Gamma(x+y)},\quad x,\;y>0,$$
we can rewrite (\ref{S5_2_Q_inv}) as
\begin{eqnarray*}
\mathcal{Q}_n(\tau)&=&\frac{1}{2\pi i}\int_{Br}\frac{e^{\tau s}}{s}\sum_{j=0}^\infty\frac{s^{-j}}{j!}\int_0^1 z^j\,(1-z)^{n+1/s}\,dz\,ds\\
&=&\frac{1}{2\pi i}\int_{Br}\frac{e^{\tau s}}{s}\int_0^1 e^{z/s}\,(1-z)^{n+1/s}\,dz\,ds\\
&=&\int_0^1 (1-z)^n\,\frac{1}{2\pi i}\int_{Br}e^{\tau s}\,\frac{e^{z/s}\,(1-z)^{1/s}}{s}\,ds\,dz.
\end{eqnarray*}
Then by the inverse Laplace transform
$$\mathcal{L}^{-1}\Big(\frac{1}{s}\,e^{W/s}\Big)=J_0(2\sqrt{\tau\,|W|}),$$
where $W=z+\log(1-z)<0$ and $J_0(\cdot)$ is the Bessel function, we obtain
$$\mathcal{Q}_n(\tau)=\int_0^1(1-z)^n\,J_0\Big(2\sqrt{\tau}\sqrt{-z-\log(1-z)}\;\Big)\,dz,$$
which leads to (\ref{S2_th5_2}). By expanding $\mathcal{Q}_n(\tau)$ for $n\to\infty$ and $\tau\to\infty$, with $\tau=O(n^2)$ we obtain
\begin{equation}\label{S5_12matching}
\mathcal{Q}_n(\tau)\sim\int_0^\infty e^{-nz}\,J_0\big(z\sqrt{2\tau}\,\big)dz=\frac{1}{\sqrt{n^2+2\tau}}=\frac{1}{\rho\,\sqrt{N^2+2T}}.
\end{equation}
Then we can easily show that $\rho^{-1}P_0(N,T)$, when expanded for $(N,T)\to(0,0)$, gives the same result as in (\ref{S5_12matching}), which verifies the matching between the long time ($T$-scale) and short time ($\tau$-scale) results.

We show how to obtain the results in Theorem 2.5 from the exact representations. First consider $n,\,t=O(\rho)$. Then from (\ref{S2_th1_nu}) and (\ref{S2_th1_nutilde}) we obtain $\nu_m\sim -m/\rho$ and $\widetilde{\nu}_m=O(\rho)$. Thus on the large time scales $t=T\rho=O(\rho)$ we have $e^{\nu_m\,t}\sim e^{-mT}$, while $e^{\widetilde{\nu}_m\,t}$ becomes exponentially small. Thus all of the terms in the first sum in (\ref{S2_th1_p}) contribute to the leading term for $p_n(t)$. We furthermore scale $n=N\rho$ and use
\begin{equation*}
C_m(\nu)\sim\frac{1}{\rho}\,\frac{m^m}{m!}\,e^{-m},\quad \left(\frac{\nu+1}{-\rho}\right)^n\sim(-\rho)^{-n}\,e^{-mN}
\end{equation*}
and
\begin{equation*}
n!\,L^{m-1-n}_n\left(\frac{\rho}{(\nu+1)^2}\right)\sim(-\rho)^n\,(1-N)^{m-1}\,e^{2mN}.
\end{equation*}
Thus from (\ref{S2_th1_p}), on the $(N,T)$ scale, we obtain
\begin{equation*}
p_n(t)\sim\rho^{-1}\,\sum_{m=1}^\infty e^{m(N-1)}\,(1-N)^{m-1}\,\frac{m^m}{m!}\,e^{-mT}.
\end{equation*}
We can show that this is equal to $\rho^{-1}P_0(N,T)$ in (\ref{S5_1_p0_series}).

Next we consider $n=O(1)$ and $t=\tau/\rho=O(\rho^{-1})$. We use the exact representation for $p_n(t)$ in Theorem 2.2 and scale $\theta=\rho s-1=O(\rho)$ ($s>0$). Then $r\sim 1/s$ and from (\ref{S2_th2_M}) we have
\begin{equation}\label{S5_2_M}
M\sim\frac{\rho^{1/s}}{s\,\Gamma(1+1/s)}.
\end{equation}
By scaling $z=w/(\rho s)=O(\rho^{-1})$ in (\ref{S2_th2_G}), $G_n$ becomes
\begin{eqnarray}\label{S5_2_Gn}
G_n&\sim&\frac{1}{(\rho s)^{n+1+1/s}}\,\int_0^1 w^n\,(1-w)^{1/s}dw\nonumber\\
&=&\frac{n!\,\Gamma(1+1/s)}{\Gamma(n+2+1/s)\,(\rho s)^{n+1+1/s}}.
\end{eqnarray}
We also find, from (\ref{S2_th2_H}), that $H_n$ is asymptotically given by
\begin{equation}\label{S5_2_Hn}
H_n\sim\int_0^\infty z^{n+1/s}\,e^{-z/s}dz=s^{n+1+1/s}\,\Gamma(n+1+1/s).
\end{equation}
Now we rewrite (\ref{S2_th2_phat}) as
\begin{equation}\label{S5_2_phat}
\widehat{p}_n(\theta)=M\,G_n\sum_{l=0}^{n-1}\frac{\rho^l}{l!}\,H_l+M\,H_n\sum_{l=n}^\infty\frac{\rho^l}{l!}\,G_l
\end{equation}
and use (\ref{S5_2_M})-(\ref{S5_2_Hn}) in (\ref{S5_2_phat}). We find that the first term in the right-hand side of (\ref{S5_2_phat}) is $O(\rho^{-2})$ and the second term is $O(\rho^{-1})$. Thus, the second term dominates the first and we have
\begin{eqnarray*}
\widehat{p}_n(\theta)&\sim&M\,H_n\sum_{l=n}^\infty\frac{\rho^l}{l!}\,G_l\\
&\sim&\frac{1}{\rho}\sum_{l=n}^\infty s^{n-l-1}\,\frac{\Gamma(n+1+1/s)}{\Gamma(l+2+1/s)}\\
&=&\frac{1}{\rho}\sum_{j=0}^\infty s^{-j-1}\,\frac{\Gamma(n+1+1/s)}{\Gamma(n+j+2+1/s)},
\end{eqnarray*}
which corresponds to $\rho^{-1}\widehat{\mathcal{Q}}_n(s)$ in (\ref{S5_2_Qhat}). Inverting the transform using $d\theta=\rho\,ds$ leads to the same result we obtained by the perturbation method.

\section{Discussion}
To summarize, we have obtained both exact and asymptotic results for the $M/M/1$-PS model with non-balking probability $b_n=1/(n+1)$. We compare our results to the standard model, where $b_n=1$. First, the spectral representation of $p_n(t)$ for the two models is very different as the standard model has a purely continuous spectrum (see also Guillemin and Boyer \cite{GU}) while the balking model has a purely discrete one.

We recently studied (see \cite{ZHmm1n}) $p_n(t)$ for the standard model asymptotically, and found that if $\rho=\lambda/\mu<1$ and $n,\;t\to\infty$ the asymptotic expansion is different according as $n/t>1-\rho$, $n/t\approx 1-\rho$, $0<n/t<1-\rho$, $n=O(t^{2/3})$, and $n=O(1)$. The scale $n=O(t^{2/3})$ is important in obtaining the tail of the unconditional density, which for the standard PS model has the form (see \cite{PO} and \cite{COH})
$$p_{_{PS}}(t)\sim\alpha_2\,t^{-5/6}\,e^{-\alpha_0\,t}\,e^{-\alpha_1\,t^{1/3}},$$
where $\alpha_0=(1-\sqrt{\rho})^2$ and $\alpha_1$ and $\alpha_2$ are constants. In contrast, for the model with balking Theorem 2.3 shows that the structure of $p_n(t)$ is different in three main sectors of the $(n,t)$ plane ($n/t>1$, $\Lambda_0<n/t<1$ and $0<n/t<\Lambda_0$), with two transition regions connecting them. For $t\to\infty$ with $0\leq n/t<\Lambda_0$ the asymptotics of $p_n(t)$ are governed by the eigenvalue with the largest real part and we obtain the purely exponential behavior in (\ref{S2_th3_5}), which leads to (\ref{S2_th6_1}) for the unconditional density $p_{_{PS}}(t)$. Thus for the model with balking the scale $n=O(t^{2/3})$ is absent.

If $\rho>1$ the standard PS model has an algebraic tail, with $p_{_{PS}}(t)\sim\alpha_3\,t^{-\rho/(\rho-1)}$, so that the mean sojourn time is finite for $\rho<2$, the second moment is finite for $\rho<3/2$, etc. Then the approximation
$$p_n(t)\sim \frac{1}{n}\,\bigg[1+(\rho-1)\,\frac{t}{n}\bigg]^{-\frac{\rho}{\rho-1}}$$
applies for $n$ and/or $t\to\infty$. This situation is similar to the model with balking in the limit $\rho\to\infty$. Here the tail will be purely exponential, but for $n$ and/or $t\to\infty$ we have the approximation in (\ref{S2_th5_1}), which is quite unlike the three sectors in Theorem 2.3.

We derived Theorem 2.5 by both a perturbation method and by using the exact representations. The former method should also be useful for general non-balking functions $b_n$, provided that we can write $\rho b_n$ in the form $\rho b_n=B(\varepsilon n)$, where $\varepsilon$ is a small parameter. Thus $\rho b_n$ is a ``slowly varying" function of $n$. For example, this would apply to $b_n=e^{-cn}$ (used by Morse \cite{morse1958}) if $c$ is small. This limit would also apply to repairman problems (or finite populations queues) where $b_n=M-n$ and $M$ is the customer population. Then $\rho b_n=\rho M (1-n/M)$ and we would assume that $M\to\infty$ (thus $\varepsilon=M^{-1}$) and $\rho\to 0$, with $\rho M=O(1)$. It is likely that the asymptotic structure of all of these models is quite different, and the perturbation method should clearly show these differences.

\end{document}